\documentclass{amsart}

\usepackage{graphicx}

\newtheorem{theorem}{Theorem}
\newtheorem{lemma}{Lemma}
\newtheorem{definition}{Definition}
\newtheorem{proposition}{Proposition}

\newcommand{\parag}[1]{\medbreak{\bf\em #1}}

\def\tilde{\widetilde}
\newcommand{\ggcd}[1]{(\!(#1)\!)}

\newcommand{\EIS}[1]{{\bf EIS:#1}}

\def\cal{\mathcal}

\def\Bbb{\mathbb}
\def\Li{\operatorname{Li}}

\def\OO{{\cal O}}
\def\ds{\displaystyle}
\def\mod{\bmod}
\def\bar{\overline}
\def\hat{\widehat}
\newcommand{\F}{{\Bbb F}}

\newcommand{\Fq}{{\F_q}}

\newcommand{\C}{{\Bbb C}}

\newcommand{\Z}{{\Bbb Z}}
\newcommand{\Q}{{\Bbb Q}}

\newcommand{\R}{{\Bbb R}}

\newcommand{\rmv}[1]{}
\def\L{\operatorname{L}}
\def\lcm{\operatorname{lcm}}
\def\asym{\operatorname{asymp}}
\def\zz{{\cal Z}}

\title[Hybrid Method]{\large\bf A Hybrid of Darboux's Method \\
and Singularity Analysis in Combinatorial Asymptotics}
\author[Flajolet, Fusy, Gourdon, Panario, and Pouyanne]
     {\rm Philippe Flajolet, Eric Fusy, Xavier Gourdon, 
     Daniel Panario, and Nicolas Pouyanne}
\address{\rm P. Flajolet, Algorithms Project, INRIA Rocquencourt,
   F-78153 Le Chesnay, France.}
\email{\tt Philippe.Flajolet{\char'100}inria.fr}
\address{\rm E. Fusy, Algorithms Project}
\email{\tt Eric.Fusy{\char'100}inria.fr}
\address{\rm X.Gourdon,
  Algorithms Project and Dassault Systems.}
\email{\tt  xgourdon{\char'100}yahoo.fr}
\address{\rm D. Panario,
School of Mathematics and Statistics, Carleton University,
   Ottawa, K1S 5B6, Canada.}
\email{\tt daniel{\char'100}math.carleton.ca}
\address{\rm N. Pouyanne,
Laboratoire de Math\'ematiques, UMR CNRS 8100,
Universit\'e de Versailles - Saint-Quentin, 78035 Versailles Cedex, France.}
\email{\tt pouyanne{\char'100}math.uvsq.fr}

\date{June 5, 2006}

\begin{document}

\begin{abstract}
A ``hybrid method'', dedicated to asymptotic coefficient extraction in
combinatorial  generating   functions, is    presented, which  combines
Darboux's method and singularity  analysis theory.  This hybrid method
applies  to  functions that remain  of  moderate  growth near the unit
circle and satisfy suitable smoothness assumptions---this, even in the
case  when the   unit   circle is  a  natural   boundary.  A   prime
application is  to    coefficients of several types  of infinite   product  generating
functions, for   which      full   asymptotic expansions   
(involving periodic fluctuations at higher orders) can      be
derived.    Examples  relative   to  permutations, trees, and
polynomials over finite fields are treated in this way.
\end{abstract}

\maketitle

\section*{Introduction}

A few enumerative problems of combinatorial  theory lead to generating
functions that are expressed  as infinite products  and admit the unit
circle  as a natural boundary.  Functions with a  fast growth near the
unit circle are usually amenable to  the saddle point method, a famous
example  being the integer  partition generating function.
We consider
here functions of moderate growth, which are  outside the scope of the
saddle point method. We do  so in the  case where neither  singularity
analysis nor Darboux's method     is directly applicable,    but the
function  to be  analysed  can  be factored  into  the product  of  an
elementary  function  with  isolated singularities  and  a sufficiently
smooth factor on  the  unit  circle.  Such decompositions   are often
attached   to infinite  products   exhibiting a regular  enough
structure and are easily    obtained by the introduction of   suitable
convergence factors.  Under  such  conditions,  we    prove   that
coefficients admit  full  asymptotic expansions  involving  powers
of logarithms and 
descending powers of  the index $n$,  as well as periodically
varying  coefficients.  Applications    are  given  to   the following
combinatorial-probabilistic problems: the enumeration of  permutations
with  distinct cycle  lengths,  the probability  that two permutations
have   the same cycle-length  profile,    the number  of  permutations
admitting an  $m$th root, the  probability  that a  polynomial over a
finite field  has factors of distinct degrees,  and the  number of forests
composed of trees of different sizes.

\smallskip
{\bf  Plan of  the      paper.}     We  start  by    recalling      in
Section~\ref{darbsa-sec}   the principles  of   two classical  methods
dedicated   to  coefficient  extraction   in  combinatorial generating
functions, namely Darboux's method and singularity analysis, which are
central to our   subsequent developments. The hybrid  method \emph{per
se}  forms  the subject of  Section~\ref{hybrid-sec},  where  our main
result,                   Theorem~\ref{hybrid2-thm},                is
established.     Section~\ref{discyc-sec}  treats   the     asymptotic
enumeration of  permutations having distinct  cycle sizes: this serves
to illustrate          in  detail      the   hybrid       method    at
work. Section~\ref{hodgepodge-sec}  discusses more  succinctly further
combinatorial problems leading to  generating functions with a natural
boundary---these    are relative     to  permutations,   forests,  and
polynomials over finite fields. A  brief perspective is offered in our
concluding section, Section~\ref{conc-sec}.

\section{Darboux's method and singularity analysis}\label{darbsa-sec}

In this section, we gather some previously known facts about Darboux's
method, singularity  analysis,  and basic   properties of analytic
functions that are central to our subsequent analyses.

\subsection{Functions of finite order}
Throughout this study,  we consider analytic functions  whose expansion
at the origin has  a finite radius  of convergence, that is, functions
with  singularities at a finite distance  from the origin. By a simple
scaling of  the independent  variable, we  may restrict   attention to
function that  are analytic in  the open unit disc $D$  but not in the
closed  unit disc $\bar  D$. What our  analysis  a priori excludes are
thus:  $(i)$~entire  functions;  $(ii)$~purely divergent  series. (For
such excluded cases,    either the  saddle  point  method   or ad  hoc
manipulations     of       divergent      series          are    often
instrumental in gaining access to coefficients~\cite{Bender74,FlSe06,Odlyzko95}.)
Furthermore we restrict attention to functions that remain
 of moderate growth near the unit circle in the following sense.

\begin{definition}\label{order-def}
A function $f(z)$ analytic in
the open unit disc~$D$ is said to be of \emph{global order}~$a\le0$  if
\[
f(z)=O((1-|z|)^a)\qquad
(|z|<1),
\]
that is, there exists an absolute constant~$M$ such that
$|f(z)|<M(1-|z|)^{a}$ for all~$z$ satisfying $|z|<1$.
\end{definition}

This definition typically excludes the partition generating function
\[
P(z)=\prod_{k=1}^\infty \frac{1}{1-z^k},
\]
which is of infinite order and to which the saddle point method
(as well as  a good deal more) is applicable~\cite{Andrews76,Ayoub63,Hardy78}.
In contrast, a function like
\[
\frac{e^z}{\sqrt{1+z}{\root 3 \of {1-z}}}\]
is of global order $a=-\frac12$, while 
\[
\exp\left(\sum_{k\ge1}\frac{z^k}{k^2}\right)
\quad{\rm or\quad}
(1-z)^{5/2}
\]
are of global order $a=0$.

We observe, though  we do not make use of the fact, that a function
$f(z)$ of global order $a\le0$ has coefficients that satisfy
$[z^n]f(z)=O(n^{-a})$.
The proof results from trivial bounds applied to Cauchy's integral form 
\begin{equation}\label{cauchy}
[z^n]f(z)=\frac{1}{2i\pi}\int_C f(z)\, \frac{dz}{z^{n+1}},
\end{equation}
upon integrating along the contour $C$: $|z|=1-n^{-1}$.
(In~\cite{BrSt97}, Braaksma and Stark present an interesting discussion
leading to refined estimates of the $O(n^{-a})$ bound.)

\subsection{Log-power functions}

What  we  address here is  the  asymptotic analysis of functions whose
local    behaviour at designated     points involves a  combination of
logarithms and powers (of possibly fractional exponent).
For the sake of notational simplicity, we write
\[
\L(z):=\log \frac{1}{1-z}.
\]
Simplifying the theory to what is needed here, we set:

\begin{definition}\label{logpow-def}
A \emph{log-power function at~$1$} is a finite sum of the form
\[
\sigma(z)=\sum_{k=1}^r c_k\left(\L(z)\right) \left(1-z\right)^{\alpha_k},\]
where $\alpha_1<\cdots<\alpha_r$ and each $c_k$ is a polynomial.
A \emph{log-power function at a finite set of points} $Z=\{\zeta_1,\ldots,\zeta_m\}$, is 
a finite sum
\[
\Sigma(z)=\sum_{j=1}^m \sigma_j\left(\frac{z}{\zeta_j}\right),
\]
where each $\sigma_j$ is a log-power function at~1.
\end{definition}
\noindent
In what follows,  we shall only  need to  consider the case  where the
$\zeta_j$ lie on the unit disc: $|\zeta_j|=1$.

\smallskip
It has been known for a long time (see, e.g., Jungen's 1931 paper, ref.~\cite{Jungen31},
and~\cite{FlOd90b,FlSe06} for wide extensions)
 that the coefficient of index~$n$  in a log-power function admits 
a full asymptotic expansion in descending powers of~$n$. 
\begin{lemma}[Coefficients of log-powers]\label{basic-lem}
The expansion of the coefficient of a log-power function is
computable by the two rules:
\begin{equation}\label{tabl}
\begin{array}{lll}
\ds [z^n](1-z)^\alpha &\sim & \ds \frac{n^{-\alpha-1}}{\Gamma(-\alpha)}
+ \frac{\alpha(\alpha+1) n^{-\alpha-2}}{\Gamma(-\alpha)} +\cdots \\
\ds [z^n](1-z)^\alpha \L(z)^k &=& \ds(-1)^k
 \frac{\partial^k}{\partial \alpha^k}
\left([z^n](1-z)^\alpha\right) \\
&\sim& \ds (-1)^k \frac{\partial^k}{\partial \alpha^k} 
\left[\frac{n^{-\alpha-1}}{\Gamma(-\alpha)}
+ \frac{\alpha(\alpha+1) n^{-\alpha-2}}{\Gamma(-\alpha)} +\cdots\right].
\end{array}
\end{equation}
\end{lemma}
The general shape of the expansion is thus
\[
\begin{array}{llll}
[z^n](1-z)^{\alpha}L(z)^k &\ds\mathop{\sim}_{n\to+\infty}&
\ds \frac{1}{\Gamma(-\alpha)}n^{-\alpha-1}(\log n)^k & (\alpha\not\in\Z_{\ge0})
\\
{}[z^n](1-z)^{r}L(z)^k &\ds\mathop{\sim}_{n\to+\infty}&
\ds (-1)^r k(r!) n^{-r-1}(\log n)^{k-1} &
(r\in\Z_{\ge0}, k\in\Z_{\ge1}).
\end{array}
\]
In the last case, the term involving $(\log n)^k$ disappears as its
coefficient is $1/\Gamma(-r)\equiv 0$.
In  essence, smaller functions  at  a singularity have  asymptotically
smaller  coefficients   and  logarithmic  factors   in  a function  are
reflected  by logarithmic terms  in  the coefficients' expansion;  for
instance,
\[
\begin{array}{lll}
\ds [z^n]\frac{\L(z)}{\sqrt{1-z}} &\sim& 
\ds \frac{\log n +\gamma+2\log 2}{\sqrt{\pi n}}-
\frac{\log n +\gamma+2\log 2}{8\sqrt{\pi n^3}}+\cdots\\
\ds[z^n](1-z)\L(z)^2 &\sim& \ds
-\frac{2}{n^2}(\log n+\gamma-1)-\frac{1}{n^3}
(2\log n+2\gamma-5)+\cdots\,.
\end{array}
\]
When supplemented by the rule
\[
[z^n]\sigma\left(\frac{z}{\zeta}\right)=\zeta^{-n}[z^n]\sigma(z),\] 
Lemma~\ref{basic-lem} 
makes it  effectively possible to determine  the asymptotic 
behaviour of coefficients of all log-power
functions.

\subsection{Smooth functions and Darboux's method}
Once the coefficients of functions in some basic scale are known,
there remains to translate error terms.
Precisely, we consider in this article functions of the form
\[
f(z)=\Sigma(z)+R(z),\] and need  conditions that enable us to estimate
the  coefficients of the  error    term $R(z)$.  Two conditions    are
classically    available:    one    based    on   smoothness    (i.e.,
differentiability)  is summarized   here, following classical  authors
(e.g.,~\cite{Olver74});   the other     based   on  growth  conditions
\emph{and} analytic continuation is discussed in the next subsection.

\begin{definition}\label{smooth-def}
Let $h(z)$ be analytic in $|z|<1$ and $s$ be a nonnegative integer.
The function $h(z)$ is said to be 
\emph{${\cal C}^s$--smooth\footnote{%
	A function $h(z)$ is said to be \emph{weakly $\mathcal{C}^s$--smooth} if it admits a continuous
	extension to the closed unit disc $|z|\le1$ and the function 
	$g(\theta):=h(e^{i\theta})$ is
	$s$~times continuously differentiable. 
	This weaker notion suffices for Lemmas~\ref{darboux-lem},
 \ref{dar-lem}< and Theorem~\ref{basic-thm}.
} on the unit disc} (or of class $\cal C^s$)
if, for all $k=0\,.\,.\,s$, its $k$th derivative
$h^{(k)}(z)$ defined for $|z|<1$ admits a  continuous extension on~$|z|\le1$. 
\end{definition}

%
For instance, a function of the form
\[
h(z)=\sum_{n\ge0} h_n z^n\qquad\hbox{with}\quad h_n=O(n^{-s-1-\delta}),
\]
for some~$\delta>0$ and $s\in\Z_{\ge0}$, is 
$\cal C^s$-smooth (both in the standard sense and in the strong sense).
Conversely, the fact that smoother functions have asymptotically smaller coefficients lies 
at the heart of Darboux's method.

\begin{lemma}[Darboux's transfer]\label{darboux-lem} 
If $h(z)$ is ${\cal C}^s$--smooth, then
\[
[z^n] h(z)=o(n^{-s}).\]
\end{lemma}
\begin{proof}
One has, by Cauchy's coefficient formula and continuity of~$h(z)$:
\[
[z^n]h(z)=\frac{1}{2\pi}\int_{-\pi}^\pi h(e^{i\theta}) e^{-ni\theta} \, d\theta.\]
When $s=0$, the statement results  directly from  the Riemann-Lebesgue
theorem~\cite[p.~109]{Rudin87}. When $s>0$, the estimate results from $s$ successive
integrations   by parts followed by the Riemann-Lebesgue
argument. 
See     Olver's
book~\cite[p.~309--310]{Olver74} for a neat discussion.
\end{proof}

\begin{definition}\label{smoothsa-def}
A function $Q(z)$ analytic  in the open unit  disc $D$ is  said to  admit a
\emph{log-power expansion of class ${\cal C}^t$} 
if there exist a finite set of points $Z=\{\zeta_1,\ldots,\zeta_m\}$
on the unit circle $|z|=1$
and a log-power function
$\Sigma(z)$ at the set of points 
$Z$ such that $Q(z)-\Sigma(z)$ is ${\cal C}^t$--smooth
on the unit circle.
\end{definition}

\begin{lemma}[Darboux's method] \label{dar-lem} If $Q(z)$ admits a log-power expansion of class~$\cal C^t$
with $\Sigma(z)$ an associated log-power function, its coefficients satisfy
\[
[z^n]Q(z)=[z^n]\Sigma(z)+o\left(n^{-t}\right).
\]
\end{lemma}

\begin{proof}
One has $Q=\Sigma+R$, with $R$ being $\cal C^t$ smooth.
The coefficients of~$R$ are estimated by Lemma~\ref{darboux-lem}.
\end{proof}


\smallskip

Consider for instance
\[
Q_1(z)=\frac{e^z}{\sqrt{1-z}}, \qquad
Q_2(z)=\frac{\sqrt{1+z}}{\sqrt{1-z}}e^z.\]
Both are of global order $-\frac12$ in the sense of Definition~\ref{order-def}. 
By making use of the analytic expansion of $e^z$ at~1, one finds
\[
Q_1(z)=\left(\frac{e}{\sqrt{1-z}}-e\sqrt{1-z}\right)+R_1(z),
\]
where $R_1(z)$, which is of the order of $(1-z)^{3/2}$ as $z\to1^-$, is
${\cal C}^1$-smooth.
The sum  of the first two terms (in parentheses) constitutes $\Sigma(z)$,
in this case with $Z=\{1\}$.
Similarly, for $Q_2(z)$, by  making use of expansions at
the elements of $Z=\{-1,+1\}$, one finds
\[
Q_2(z)=\left(\frac{e\sqrt{2}}{\sqrt{1-z}}-\frac{5e}4 \sqrt{2}\sqrt{1-z}
+\frac{1}{e\sqrt{2}}\sqrt{1+z}\right)+R_2(z),
\]
where  $R_2(z)$ is also  ${\cal C}^1$--smooth. 
Accordingly, we find:
\begin{equation}\label{p1np2n}
[z^n]Q_1(z)=e\frac{1}{\sqrt{\pi n}}+o\left(\frac{1}{n}\right),
\qquad
[z^n]Q_2(z)=\frac{e\sqrt{2}}{\sqrt{\pi  n}}+o\left(\frac{1}{n}\right).
\end{equation}  The  next    term in  the     asymptotic expansion of
$[z^n]Q_2$   involves   a linear  combination   of   $n^{-3/2}$  and $(-1)^n
n^{-3/2}$, where the latter term reflects the singularity  at $z=-1$.  Such
calculations are typical of what we shall encounter later.

%

\subsection{Singularity analysis}

What we refer to as singularity  analysis is a technology developed by
Flajolet and Odlyzko~\cite{FlOd90b,Odlyzko95}, with further additions
to be found in~\cite{FiFlKa05,Flajolet99,FlSe06}.
It applies to a function with a finite number of 
singularities on the boundary of its disc of convergence. 
Our description closely follows Chapter~VI of the latest edition
of \emph{Analytic Combinatorics}~\cite{FlSe06}.

Singularity analysis  theory adds to  Lemma~\ref{basic-lem}
the theorem that, under  conditions  of analytic continuation, $O$-  and
$o$-error terms can  be similarly transferred to coefficients.
Define a \emph{$\Delta$-domain} associated to two parameters $R>1$ (the radius) and 
$\phi\in(0,\frac{\pi}{2})$ (the angle) by
\[
\Delta(R,\phi):=\left\{ z ~\bigm|~|z|<R,~ \phi<\arg(z-1)<2\pi-\phi, ~z\not=1\right\}
\]
where $\arg(w)$ denotes the argument of $w$ taken here in the interval $[0,2\pi[$.
By definition a $\Delta$-domain properly  contains the unit disc, since $\phi<\frac{\pi}{2}$.
(Details of the values of $R,\phi$ are immaterial as long as $R>1$ and $\phi<\frac\pi2$.)

The following definition is in a way the counterpart of smoothness 
(Definition~\ref{smoothsa-def}) for singularity analysis
of functions with isolated singularities.

\begin{definition} \label{ordersa-def}
Let $h(z)$ be analytic in $|z|<1$ and have isolated singularities
on the unit circle at $Z=\{\zeta_1,\ldots,\zeta_m\}$.
Let $t$ be a real number.
The function  $h(z)$ is said to admit a \emph{log-power expansion of type $\OO^t$}
(relative to~$Z$) if the following two conditions are satisfied:
\begin{itemize}
\item[---] The function $h(z)$ is analytically continuable to
an indented domain $\mathfrak{D}=\bigcap_{j=1}^m (\zeta_j\cdot \Delta)$,
with $\Delta$ some  $\Delta$-domain. 
\item[---] There exists a log-power function $\Sigma(z):=\sum_{j=1}^m \sigma_j(z/\zeta_j)$
such that, for each $\zeta_j\in Z$, one has
\begin{equation}\label{sacond}
h(z)-\sigma_j(z/\zeta_j)=O\left((z-\zeta_j)^t\right),
\end{equation}
as $z\to\zeta_j$ in $(\zeta_j\cdot\Delta)$.
\end{itemize}
\end{definition}
Observe that  $\Sigma(z)$ is \emph{a priori}  uniquely determined only up
to $O((z-\zeta_j)^t)$ terms. The  minimal function (with respect to the
number    of  monomials)  satisfying~\eqref{sacond}   is   called  the
\emph{singular part} of $h(z)$ (up to $\OO^t$  terms).

A basic result of singularity analysis theory
enables us to extract coefficients of functions that admit of such expansions.

\begin{lemma}[Singularity analysis method] \label{sa-lem}
Let $Z=\{\zeta_1,\ldots,\zeta_m\}$ be a finite set of points on the unit circle,
and let $P(z)$ be a function that admits
a log-power expansion of type $\cal O^t$ relative to $Z$,
with singular part $\Sigma(z)$.
Then, the coefficients of~$h$ satisfy
\begin{equation}\label{sares}
[z^n]P(z)=[z^n]\Sigma(z)+O\left(n^{-t-1}\right).
\end{equation}
\end{lemma}
\noindent
\begin{proof}
The proof of Lemma~\ref{sa-lem} starts from Cauchy's integral formula~\eqref{cauchy}
and makes use of the contour $C$ that lies
at distance $\frac1n$ of the boundary of the analyticity domain,
$\mathfrak{D}=\bigcap_{j=1}^m \left(\zeta_j\cdot \Delta\right)$. 
See~\cite{FlOd90b,FlSe06} for details.
\end{proof}

\subsection{Polylogarithms} 
For future reference (see especially Section~\ref{discyc-sec}), we gather here facts
relative to the \emph{polylogarithm} function $\Li_\nu(z)$,
which is defined for any $\nu\in\C$ by
\begin{equation}\label{polydef}
\Li_\nu(z):=\sum_{n=1}^\infty \frac{z^n}{n^\nu}.
\end{equation}
One has in particular
\[
\Li_0(z)=\frac{z}{1-z},\qquad
\Li_1(z)=\log\frac{1}{1-z}\equiv  \L(z).
\]
In the most basic applications, one encounters polylogarithms of integer index,
but in this paper (see the example of dissimilar forests in Section~\ref{hodgepodge-sec}),
the more general case of a real index~$\nu$ is also needed.

\begin{lemma}[Singularities of polylogarithms]
For any index~$\nu\in\C$, the polylogarithm~$\Li_\nu(z)$ is analytically
continuable to the slit plane $\C\setminus\R_{\ge1}$.
If~$\nu=m\in\Z_{\ge1}$, the singular expansion of $\Li_m(z)$ 
near the singularity $z=1$ is given by
\begin{equation}\label{polyexp}\renewcommand{\arraycolsep}{3truept}
\left\{
\begin{array}{rll}
\Li_{m}(z)
&=&\ds \frac{(-1)^m}{(m-1)!}\tau^{m-1}(\log \tau-H_{m-1})+
\!\!\sum_{j\ge0, j\not=m-1}
\frac{(-1)^j}{j!}\zeta(m-j)\tau^j
\\
\tau& :=& \ds -\log z \quad = \quad \sum_{\ell=1}^\infty \frac{(1-z)^\ell}{\ell}.
\end{array}\right.
\end{equation}
For $\nu$ not an integer, the singular expansion of $\Li_\nu(z)$ is
\begin{equation}\label{polyexp2}
\Li_\nu(z)\sim \Gamma(1-\nu)\tau^{\nu-1}+\sum_{j\ge0}\frac{(-1)^j}{j!}\zeta(\nu-j)\tau^j.
\end{equation}
\end{lemma}
The representations are given as
a \emph{composition} of two explicit series.
The expansions involve both the harmonic number $H_m$ and the Riemann zeta function $\zeta(s)$ defined by
\[
H_m=1+\frac12+\frac13+\cdots+\frac{1}{m},\qquad
\zeta(s)=\frac{1}{1^s}+\frac{1}{2^s}+\frac{1}{3^s}+\cdots\,
\]
($\zeta(s)$, originally defined in the half-plan $\Re (s)>1$, is
analytically continuable to $\C\setminus\{ 1\}$ by virtue of its classical functional
equation).
\begin{proof}

First in the case of an integer index~$m\in\Z_{\ge2}$,
since $\Li_m(z)$ is an iterated integral of $\Li_1(z)$,
it is analytically continuable to the complex plane slit along the ray $[1,+\infty[$. By this device, its
expansion at the singularity $z=1$ can be determined, resulting in~\eqref{polyexp}.
(The representation in~\eqref{polyexp} is in fact exact and not
merely asymptotic. It has been 
obtained by Zagier and Cohen in~\cite[p.~387]{Lewin91},
and is known to the symbolic manipulation system {\sc Maple}.)

For $\nu$ not an integer,  analytic continuation derives from
a Lindel\"of integral representation discussed by Ford in~\cite{Ford60}.
The singular expansion, valid in the slit plane, was 
established in~\cite{Flajolet99} to which we refer for details.
\end{proof}

In the sequel, we also make use of smoothness properties of polylogarithms.
Clearly, $\Li_k(z)$ is $\cal C^{k-2}$--smooth in the sense of Definition~\ref{smooth-def}.
A simple computation of coefficients shows that any sum
\[
S_{k}(z)=\sum_{\ell\ge k} r(\ell) \left[\Li_\ell(z^\ell)-\Li_{\ell}(1)\right]
\]
with $r(x)$ polynomially bounded in~$x$, is $\cal C^{k-2}$--smooth.  Many similar
sums     are    encountered   later,   starting     with    those   in
Equations~\eqref{maindcyc} and~\eqref{maindcyc2}.

\section{The hybrid method} \label{hybrid-sec}

The heart  of the matter  is  the treatment of functions  analytic in the
open unit disc that can, at least partially, be ``de-singularized'' by
means of log-power functions.

\subsection{Basic technology}
Our first theorem, which essentially relies on the Darboux technology,
serves as a  stepping stone towards the proof  of  our main statement,
Theorem~\ref{hybrid2-thm} below.

\begin{theorem}\label{basic-thm}
Let $f(z)$ be analytic in the open unit disc $D$, of global
order   $a\le  0$, and such that it admits a  factorization
$f=P\cdot Q$,   with  $P,Q$  analytic in~$D$.  Assume    the following
conditions on $P$ and $Q$, relative  to a finite  set of points
$Z=\{\zeta_1,\ldots,\zeta_m\}$ on the unit circle:
\begin{itemize}
\item[${\bf C}_1$:] The ``Darboux factor''  $Q(z)$ is ${\cal C}^s$--smooth
on the unit circle ($s\in\Z _{\geq 0}$).
\item[${\bf C}_2$:] The ``singular factor'' $P(z)$ admits,
for some nonnegative integer $t$,
a log-power expansion relative to~$Z$, $P=\tilde P+R$ (with $\tilde P$ the log-power
function and $R$ the smooth term), that is of class
${\cal C}^t$. 
\end{itemize}
Assume also the inequality (with~$\lfloor x\rfloor$ the integer part function):
\begin{itemize}
\item[${\bf C_3}$:] 
$t\ge u_0$, where
\begin{equation}\label{defu0}
u_0:=\left\lfloor \frac{s+\lfloor a\rfloor}{2}\right\rfloor.
\end{equation}
\end{itemize}
Let $c_0=\left\lfloor \frac{s-\lfloor a\rfloor}{2}\right\rfloor$.
If $H$ denotes the Hermite interpolation polynomial\footnote{
	Hermite  interpolation  extends the  usual process of Lagrange
	interpolation,  by allowing   for  higher  contact   between a
	function and its interpolating  polynomial at a designated set
	of  points.  A  lucid construction  is  found in
	Hildebrand's treatise~\cite[\S8.2]{Hildebrand74}.
} such that all its
derivatives of order
$0,\dots ,c_0-1$ 
coincide with those of $Q$ at each of the points $\zeta _1,\dots ,\zeta _m$,
one has
\begin{equation}\label{conclthm1}
[z^n]f(z)=[z^n]\left(\tilde P(z)\cdot H(z)\right)+o(n^{-u_0}).
\end{equation}
\end{theorem}
\noindent
Since $\tilde P(z)\cdot H(z)$  is itself a log-power function, the asymptotic
form     of   its     coefficients is   explicitly   provided    by
Lemma~\ref{basic-lem}.

\begin{proof}
Let~$c\le s$  be a positive integer whose precise value
will be adjusted at the end of the proof.  First, we decompose $Q$ as
\[
Q=\bar Q+S,
\]
where $\bar Q$  is the polynomial of
minimal degree such that all its  derivatives of order $0,\ldots,c-1$ at
each  of  the  points   $\zeta_1,\ldots,\zeta_m$ coincide   with those
of~$Q$:
\begin{equation}\label{herm1}
\left.\frac{\partial^i}{\partial z^i} \bar Q(z)\right|_{z=\zeta_j}
=
\left.\frac{\partial^i}{\partial z^i}  Q(z)\right|_{z=\zeta_j},
\qquad 0\le i< c,\quad
1\le j\le m.
\end{equation}
The   classical process  of Hermite  interpolation~\cite{Hildebrand74}
produces such a polynomial, whose degree is at most $cm-1$.  Since $\bar Q$ is
${\cal  C}^\infty$--smooth,   the    quantity  $S=Q-\bar Q$     is ${\cal
C}^s$--smooth. This function~$S$ is also ``flat'',   in the sense that
it has a contact of high order with~0 at each of the points~$\zeta_j$.

We now operate with the decomposition
\begin{equation}\label{majdecomp}
f=\tilde P \cdot \bar Q + \tilde P \cdot S + R\cdot \bar Q + R\cdot S,
\end{equation}
and proceed to examine the coefficient of $z^n$ in each term.

\smallskip
---~\emph{The product $\tilde P \cdot \bar Q$.} 
Since  $\tilde P$ is a log-power  function and $\bar  Q$ a polynomial,
the    coefficient    of   $z^n$  in    the    product   admits,    by
Lemma~\ref{basic-lem}, a complete descending  expansion with  terms in
the scale $\{n^{-\beta}(\log n)^k\}$, which we write concisely as
\begin{equation}\label{red1}
[z^n]\tilde P \cdot \bar Q\in \left\{
n^{-\beta}(\log n)^k~\bigm|~k\in\Z_{\ge0},~\beta\in \R\right\}.
\end{equation} 

\smallskip
---~\emph{The product $\tilde P \cdot S$.}
This is where the Hermite interpolation polynomial $\bar Q$ plays its part.
From the construction of $\bar Q$, there results that 
$S=Q-\bar Q$ has all its derivatives of order $0,\ldots,c-1$ vanishing at each of the points 
$\zeta_1,\ldots,\zeta_m$.
This guarantees the existence of a factorization
\[
S(z)\equiv Q(z)-\bar Q(z) = \kappa(z)\prod_{j=1}^m (z-\zeta_j)^c,
\]
where $\kappa(z)$ is now  ${\cal C}^{s-c}$--smooth (division  decreases the
degree of smoothness). Then, in the factorization
\[
\tilde P \cdot  S =\left(\tilde P \cdot \prod_{j=1}^m (z-\zeta_j)^c \right)\cdot \kappa(z),
\]
the quantity $\tilde P$ is, near a~$\zeta_j$,
of order at most $O(z-\zeta_j)^a$ (with $a$ the global order of $f$).
Thus, $\tilde P S/\kappa$ is at least ${\cal C}^v$--smooth,
with $v:=\lfloor c+a\rfloor$. Since ${\cal C}^p\cdot{\cal C}^q\subset {\cal C}^{\min(p,q)}$,
Darboux's method (Lemma~\ref{dar-lem}) yields
\begin{equation}\label{red2}
[z^n]\widetilde P \cdot S =o\left(n^{-u(c)}\right), \qquad
u(c):=\min(\lfloor c+a\rfloor,s-c).
\end{equation}

\smallskip
---~\emph{The product $R \cdot \bar Q$.}
This   quantity  is of   class  ${\cal  C}^t$  by  elementary  product
rules. Hence, by Darboux's method,
\begin{equation}\label{red3}
[z^n] R\cdot \bar Q=o\left(n^{-t}\right).
\end{equation}

\smallskip
---~\emph{The product $R \cdot  S$.}
This product is of class ${\cal C}^{\min(s,t)}$ and, by Darboux's method again,
\begin{equation}\label{red4}
[z^n] R\cdot S=o\left(n^{-\min(s,t)}\right).
\end{equation}

\smallskip

It now only remains to collect the effect of the various error terms
of~\eqref{red2}, \eqref{red3}, and~\eqref{red4} in the
decomposition~\eqref{majdecomp}:
\[
[z^n]f = \left([z^n] \tilde P \cdot \bar Q \right) + 
o(n^{-u(c)})+o(n^{-t})+o(n^{-\min(s,t)}).
\]
Given the condition $t\ge u_0$ in ${\bf C_3}$, the last two terms are
$o(n^{-u_0})$. 
A choice, which  maximizes $u(c)$ (as defined in~\eqref{red2}) 
and suffices for our purposes, is
\begin{equation}\label{optc}
c_0=\left\lfloor \frac{s-\lfloor a\rfloor}{2}\right\rfloor
\qquad\hbox{corresponding to}\qquad
u(c_0)=\left\lfloor \frac{s+\lfloor a\rfloor}{2}\right\rfloor=u_0.
\end{equation}
The statement then results from  the choice of $c=c_0$, as well as
$u_0=u(c_0)$ and $H(z):=\bar Q(z)$, the corresponding Hermite interpolation
polynomial.
\end{proof}

\subsection{Hybridization}

Theorem~\ref{basic-thm} is largely to   be  regarded as an   existence
result: due to  the factorization and the  presence  of a Hermite
interpolation polynomial, it  is  not   well suited  for   effectively
deriving  asymptotic expansions.   In this subsection,  we develop the
hybrid method  \emph{per se},  which  makes  it  possible  to  operate
directly  with  a small number of   radial expansions  of the function
whose coefficients are to be estimated.

\begin{definition} Let $f(z)$ be analytic in the open unit disc.
For $\zeta$ a point on the unit circle, we define the 
\emph{radial expansion} of $f$ at~$\zeta$ with order $t\in\R$ as the
smallest (in terms of the number of monomials)
log-power function $\sigma(z)$ at~$\zeta$,
provided it exists, such that
\[
f(z)=\sigma(z)+O\left((z-\zeta)^t\right),
\]
when $z=(1-x)\zeta$ and $x$ tends to~$0^+$.
The quantity $\sigma(z)$ is written 
\[
\asym(f(z),\zeta,t).
\] 
\end{definition}
The interest of radial expansions is to a large extent a computational
one,  as these are  often accessible via  common methods of asymptotic
analysis while various series  rearrangements   from within the   unit
circle  are granted  by   analyticity.    In contrast,   the   task of
estimating directly a function $f(z)$ as $z\to\zeta$ \emph{on the unit
circle} may be technically more demanding.  Our main theorem is accordingly
expressed in terms of such radial expansions  and, after the necessary
conditions on the generating function  have been verified, it provides
an
\emph{algorithm} (Equation~\eqref{algo})  for the determination of the
asymptotic form of coefficients.

\begin{theorem}[Hybrid method]\label{hybrid2-thm}
Let $f(z)$  be  analytic  in the open   unit disc  $D$,
of finite global order $a\le 0$, and such that it
admits  a
factorization  $f=P\cdot  Q$, with $P,Q$ analytic in~$D$.
Assume    the following  conditions  on
$P$ and $Q$, relative to a finite set of points $Z=\{\zeta_1,\ldots,\zeta_m\}$
on the unit circle:
\begin{itemize}
\item[${\bf D}_1$:] 
The ``Darboux  factor''  $Q(z)$ is  ${\cal  C}^{s}$--smooth  on the unit
circle ($s\in\Z _{\geq 0}$).
\item[${\bf D}_2$:] 
The   ``singular  factor''  $P(z)$ is analytically continuable to an indented
domain of the form
$\mathfrak{D}=\bigcap_{j=1}^m \left(\zeta_j \cdot \Delta\right)$. 
For some
non-negative real number $t_0$, it admits, at any
$\zeta_j\in Z$, an asymptotic expansion
of the form 
\[
P(z)=\sigma_j(z/\zeta_j)+O\left((z-\zeta_j)^{t_0}\right) \qquad
(z\to\zeta_j,~z\in\mathfrak{D}),
\]
where  $\sigma_j(z)$ is a log-power function at~$1$.
\end{itemize}
Assume also the inequality:
\begin{itemize}
\item[$\bf D_3$:] 
$t_0>\lfloor \frac{s+\lfloor a\rfloor}2\rfloor$.
\end{itemize} 
Then $f$ admits a radial expansion at any $\zeta _j\in Z$ with order
$u_0=\left\lfloor \frac{s+\lfloor a\rfloor}{2}\right\rfloor$.
The coefficients of~$f(z)$ satisfy:
\begin{equation}\label{algo}
\begin{array}{l}
[z^n]f(z) \quad = \quad \ds [z^n]{\bf A}(z)+o\left(n^{-u_0}\right),
\\
\hbox{where} \qquad \ds {\bf  A}(z):=\sum_{j=1}^m      \asym(f(z),\zeta_j,u_0).
\end{array}
\end{equation}
\end{theorem}

\begin{proof}
Let us denote by
\[
\Sigma(z)=\sum_{j=1}^m \sigma_j(z/\zeta_j),
\]
the sum  of the  singular  parts of  $P$  at the points  of $\zeta_j$.
The difference $R:=P-\Sigma$ is ${\cal C}^{t}$-smooth 
for any integer~$t$ satisfying~$t<t_0$ (in particular, we can choose~$t=u_0$,
this by assumption~${\bf D_2}$.
The singular factor $P$ has thus been re-expressed as the sum of a singular part~$\Sigma$
and a smooth part $R$.
The  conditions of  Theorem~\ref{basic-thm} are then precisely satisfied
by the product $PQ$, the inequality ${\bf D_3}$ implying condition~$\bf C_3$,
so that one has by~\eqref{conclthm1}
\begin{equation}\label{above}
[z^n]f(z)=[z^n]\Sigma(z)H(z)+o(n^{-u_0}),
\end{equation}
where $H$ is the Hermite polynomial associated with $Q$ that is described in the proof of
Theorem~\ref{basic-thm} and $u_0$ is given by~\eqref{defu0}.

In order to complete
the proof,  there   remains  to  verify  that, in the  coefficient
extraction  process of~\eqref{above} above,
the quantity $\Sigma H$   can  be replaced by  $\mathbf{A}(z)$.

We have
\begin{equation}\label{interim1}
[z^n]\Sigma(z)H(z) = \sum_j [z^n]\sigma_j(z/\zeta_j)H(z).
\end{equation}
Now, near each $\zeta_j$, we have (with~$c_0=\left\lfloor \frac{s-\lfloor a\rfloor}{2}\right\rfloor$ according to~\eqref{optc})
\begin{equation}\label{bounds}
\begin{array}{lll}
\sigma_j(z/\zeta_j)&=&P(z)+O((z-\zeta_j)^{t_0})\\
H(z)&=&Q(z)+O((z-\zeta_j)^{c_0})\\
P(z)&=&O((z-\zeta_j)^a),
\end{array}
\end{equation}
respectively by assumption~$\bf D_2$, by the high order contact of~$H$ with~$Q$
due to the Hermite interpolation construction, and by the global order property of~$f(z)$.
There results from Equation~\eqref{bounds}, condition~$\bf D_3$,
and the value of~$c_0$ in~\eqref{optc} that
\[
\sigma_j(z/\zeta_j)H(z)=\asym(P(z)Q(z),\zeta_j,u_0)+O((z-\zeta_j)^{u_0}),\]
The proof, given~\eqref{above} and \eqref{interim1}, is now complete.
\end{proof}
%
%
%
%
%
%

%

Thanks to Theorem~\ref{hybrid2-thm}, in order to analyse the coefficients
of a function~$f$, the following two steps are sufficient.
\begin{itemize}
\item[$(i)$] Establish the \emph{existence} of a proper factorization~$f=P\cdot Q$.
Usually, a crude analysis is sufficient for this purpose.
\item[$(ii)$] Analyse \emph{separately} the asymptotic character of 
$f(z)$ as $z$ tends radially to a few distinguished points, those of~$Z$.
\end{itemize}
As asserted by Theorem~\ref{hybrid2-thm},  it then becomes    possible  to
proceed with  the analysis of  the coefficients $[z^n]f(z)$ \emph{as though} 
the function~$f$  satisfied
the conditions of singularity analysis (whereas in general $f(z)$ admits
the unit circle as a natural boundary).

\smallskip

Manstavi{\v c}ius~\cite{Manstavicius02} develops an alternative approach that requires conditions on generating functions in 
the disc of convergence, but only weak smoothness on the circumference. His results are
however not clearly adapted to deriving symptotic expansions beyond the main terms.

\section{Permutations with distinct cycle sizes}\label{discyc-sec}

%
%

The function
\[
f(z):=\prod_{k=1}^\infty \left(1+\frac{z^k}{k}\right),
\]
has been studied by Greene and  Knuth~\cite{GrKn81}, in relation to
a problem relative to factorization  of polynomials over finite fields
that we treat later.   As is readily  recognized from first principles
of  combinatorial  analysis~\cite{FlSe06,GoJa83,Stanley98,Wilf94}, the
coefficient  $[z^n]f(z)$  represents  \emph{the   probability that, in
a random  permutation of size    $n$, all cycle lengths are   distinct.}
One has
\[
f(z)=1+z+\frac{z^2}{2!}+5\frac{z^3}{3!}+14\frac{z^4}{4!}+74\frac{z^5}{5!}+\cdots,
\]
and the coefficients constitute the sequence~\EIS{A007838}\footnote{
	We shall use the notation~\EIS{xxxxxx} to represent a sequence indexed 
	in the \emph{Encyclopedia of Integer Sequences}~\cite{Sloane06}.}.
In~\cite[\S4.1.6]{GrKn81},    the  authors  devote  some  seven  pages
(pp.~52--58) to the derivation of the estimate ($\gamma$ is Euler's constant)
\begin{equation}\label{grkn} 
f_n:=[z^n]f(z)
=e^{-\gamma}+\frac{e^{-\gamma}}{n}+O\left(\frac{\log n}{n^2}\right),
\end{equation}
starting with a Tauberian argument and repeatedly using bootstrapping.
In our treatment below, we recycle some  of their calculations, though
our asymptotic technology is fundamentally different.

\smallskip
\emph{Global order.}
The first task in our perspective is to determine the global order
of $f(z)$.
The following chain of calculations,
\begin{equation}\label{chain}
\begin{array}{lll}
f(z)&=&\ds \prod_{k=1}^\infty e^{z^k/k}
\prod_{k=1}^\infty \left(1+\frac{z^k}{k}\right)e^{-z^k/k}\\
&=& \ds \frac{1}{1-z}\exp\left(
\sum_{k=1}^\infty \log\left(1+\frac{z^k}{k}\right)-\frac{z^k}{k}\right)\\
&=&\ds \frac{1}{1-z}\exp\left(-\frac12\sum_{k\ge1}\frac{z^{2k}}{k^2}
+\frac13\sum_{k\ge1}\frac{z^{3k}}{k^3}-\cdots\right),
\end{array}
\end{equation}
shows $f(z)$ to be of global order~$-1$. 
It is based on the usual introduction of convergence factors, the exp--log
transformation ($X\equiv\exp(\log X)$), and finally the logarithmic expansion.

Note that this preliminary determination of global order only gives the
useless bound $f_n=O(n )$.
Actually, from the infinite product expression of the Gamma
function~\cite{WhWa27}
(or from a direct calculation, as in~\cite{GrKn81}), there results that 
\begin{equation}\label{euler}
e^{-\gamma}=\prod_{k\ge1}\left(1+\frac{1}{k}\right)e^{-1/k},
\end{equation}
hence, from the second line of~(\ref{chain}), 
\[
f(z)\mathop{\sim}_{z\to1^-} \frac{e^{-\gamma}}{1-z},
\]
which is compatible with \eqref{grkn}, but far from sufficient to imply it.

\smallskip
\emph{The hybrid method.}
Given the last line of~\eqref{chain}, 
which is re-expressed in terms of polylogarithms as
\begin{equation}\label{maindcyc}
f(z)=
\frac{1}{1-z}\exp\left(-\frac12\Li_2(z^2)+\frac{1}{3}\Li_3(z^3)-\cdots\right),
\end{equation}
the right factorization of $f(z)$ is obtained transparently. Define
\begin{equation}\label{maindcyc2}
U(z):=\sum_{1\le \ell\le s+1} \frac{(-1)^{\ell-1}}{\ell}\Li_\ell(z^\ell),
\qquad
V(z):=\sum_{s+2\le \ell} \frac{(-1)^{\ell-1}}{\ell}\Li_\ell(z^\ell),
\end{equation}
so that 
\begin{equation}\label{factotum}
f(z)=e^{U(z)}\cdot e^{V(z)}.
\end{equation}
Clearly $V(z)$ is ${\cal C}^s$--smooth, and so is $e^{V(z)}$ given usual rules
of differentiation. Thus $Q:=e^V$ is our Darboux factor. The first factor $P:=e^U$ 
satisfies the condition of Theorem~\ref{hybrid2-thm}:
it is the singular factor  and it can be expanded to any order~$t$
of smallness. Consequently, the hybrid method is applicable and  can provide an asymptotic 
expansion of $[z^n]f(z)$ to any predetermined degree of accuracy.

%

\smallskip
\emph{The nature of the full expansion.} 
Given the existence of factorizations of type~\eqref{factotum} with an
arbitrary degree of smoothness (for  $V$) and smallness (for $U$),  it
is possible  to   organize the  calculations   as  follows: take   the
primitive roots of  unity   in sequence, for   orders  $1,2,3,\ldots$.
Given such a  root~$\eta$  of order~$\ell$, each  radial restriction  admits a
full asymptotic expansion in descending powers of $(1-z/\eta)$ tempered
by  polynomials  in   $\log(1-z/\eta)$.    Such an  expansion  can   be
translated  formally into a
full expansion in powers of $n^{-1}$ tempered  by polynomials in $\log
n$ and multiplied  by $\eta^{-n}$. All the  terms collected in  this way
are bound to occur in the asymptotic expansion of $f_n=[z^n]f(z)$.

For the sequel, it proves convenient first to adjust 
the expansion~\eqref{chain}, by taking out the $(1+z)$ factor. We find,
by the same techniques as in~\eqref{chain}
\begin{equation}\label{chain2}
f(z)=e^{-z}\frac{1+z}{1-z} \exp\left(\sum_{m\ge2}
\frac{(-1)^{m-1}}{m}\left[\Li_m(z^m)-z^m\right]\right),
\end{equation}
and start with the expansion as $z\to1$, then consider 
in turn $z=-1$,
$z=\omega,\omega^2$ (with $\omega=\exp(2i\pi/3)$),
and finally $z=\eta$, a primitive $\ell$th root of unity.

\smallskip

\emph{The expansion at~$z=1$.}
Calculations simplify a bit if we set
\[
z=e^{-\tau}, \qquad \tau=-\log z,
\]
as in~\eqref{polyexp}. 
By summing the  singular expansions of polylogarithms~\eqref{polyexp},
one arrives at an asymptotic expansion as $\tau\to0^+$ of the form:
\begin{equation}\label{main1}
f(e^{-\tau})=e^{-e^{-\tau}}\frac{1+e^{-\tau}}{1-e^{-\tau}}
\exp\left(-\alpha(\tau)\log\tau+\beta(\tau)+\delta_1(\tau)+\delta_2(\tau)-\epsilon(\tau)\right).
\end{equation}
There, the first two terms in the exponential correspond to summing the 
special terms in the singular expansions of polylogarithms~\eqref{polyexp}:
\[
\alpha(\tau):=\sum_{m\ge1} \frac{(m+1)^{m-1}\tau^m}{m!},
\qquad
\beta(\tau)=\sum_{m\ge1} \frac{(m+1)^{m-1}\tau^m}{m!}(H_m-\log(m+1)).\]
The last three terms 
inside the exponential of~(\ref{main1})
arise from summation over values of $m\ge2$ of the regular part 
of~\eqref{polyexp}, namely,
\[
\sum_{m\ge2}\frac{(-1)^{m-1}}{m}
\sum_{j\ge0,~j\not=m-1}\frac{(-1)^j}{j!}\zeta(m-j)(m\tau)^j,
\]
upon distinguishing between the three cases:
$j<m-1$ (giving rise to $\delta_1(\tau)$),
 $j=m-1$ (giving $\epsilon(\tau)$), and $j>m-1$ (giving  $\delta_2(\tau)$),
and exchanging the order of summations.
The calculation of $\epsilon(\tau)$ and $\delta_1(\tau)$ is immediate.
First
\begin{equation}\label{del1}
\delta_1(\tau)=\sum_{j\ge 2}\delta_{1,j}\frac{(-\tau)^j}{j!},
\qquad
\delta_{1,j}=\sum_{m=2}^j(-1)^{m-1}m^{j-1}(\zeta(m-j)-1),
\end{equation}
and since $\zeta$ values at negative integers are rational numbers,
the expansion of $\delta_1(\tau)$ involves only rational coefficients.
Next, one finds
\[
\epsilon(\tau)=\sum_{m\ge1} (m+1)^{m-1}\frac{\tau^m}{m!},
\]
a variant of the Lambert and Cayley functions. Finally, 
the function $\delta_2(\tau)$ has coefficients \emph{a priori} 
given by sums like in~\eqref{del1}, but with the summation extending
to $m\ge j+2$:
\[
\delta_2(\tau)=\sum_{j\ge0} \delta_{2,j}\frac{(-\tau)^j}{j!},
\qquad
\delta_{2,j}=\sum_{m\ge j+2} (-1)^{m-1}m^{j-1}(\zeta(m-j)-1).
\]
Each infinite sum in the expansion of $\delta_2$
is expressible in finite form: it suffices to 
start from the known expansion of $\psi(1+s)$ at $s=0$, which gives
($\psi(s)$ is the logarithmic derivative of the Gamma function)
\begin{equation}\label{psiexp}
\psi(1+s)+\gamma-\frac{s}{1+s}=(\zeta(2)-1)s-(\zeta(3)-1)s^2+\cdots\,,
\end{equation}
and differentiate an arbitrary number of times with respect to~$s$, then 
finally set $s=1$.
One finds for instance, in this way,
\[
\delta_{2,0}=-\gamma-\log2+1,\quad
\delta_{2,1}=\frac12,\quad
\delta_{2,2}=2\zeta(3)-\frac{\pi^2}{2}+\frac32.\]

Equation~\eqref{main1}   provides   a complete   algorithm  for
expanding $f(z)$ as $z\to1^{-}$. 
The first few terms found are
\begin{small}
\[
\begin{array}{rcl}
f(z)&=&
\ds e^{-\gamma}\bigg(\frac{1}{1-z}-\log(1-z)-\log 2
+\frac{1}{2}(1-z)\log^2(1-z)\\
&&\hskip 50pt\ds +(\log 2-2)(1-z)\log(1-z)+O(1-z)bigg).
\end{array}
\]
\end{small}%
Then, an application of Theorem~\ref{hybrid2-thm} yields the terms in the asymptotic expansion 
of $f_n$ arising from the singularity $z=1$:
\begin{equation}\label{at1}
\hat f_n^{[1]}=e^{-\gamma}+\frac{e^{-\gamma}}{n}+\frac{e^{-\gamma}}{n^2}
(-\log n+c_{2,0})+\frac{e^{-\gamma}}{n^3}(\log^2 n+c_{3,1}\log n+c_{3,0})+\cdots\,
\end{equation}
where
\begin{equation}\label{c30eq}
\begin{array}{l}
c_{2,0}=-1-\gamma+\log2,
\qquad
c_{3,1}=4+2\gamma-2\log2,
\\
c_{3,0}= 
1+4\gamma-\log2-3\log3+\log^2 2-\frac{\pi^2}{3}+\gamma^2-2\gamma\log2.
\end{array}
\end{equation}
From preceding considerations, the coefficients all lie in the ring generated by $\log2,\log3,\ldots$, $\gamma$, $\pi$, and $\zeta(3),\zeta(5),\ldots$.

\begin{figure}

\begin{center}
\hbox{\hspace*{1.5truecm}\includegraphics[width=50truemm]{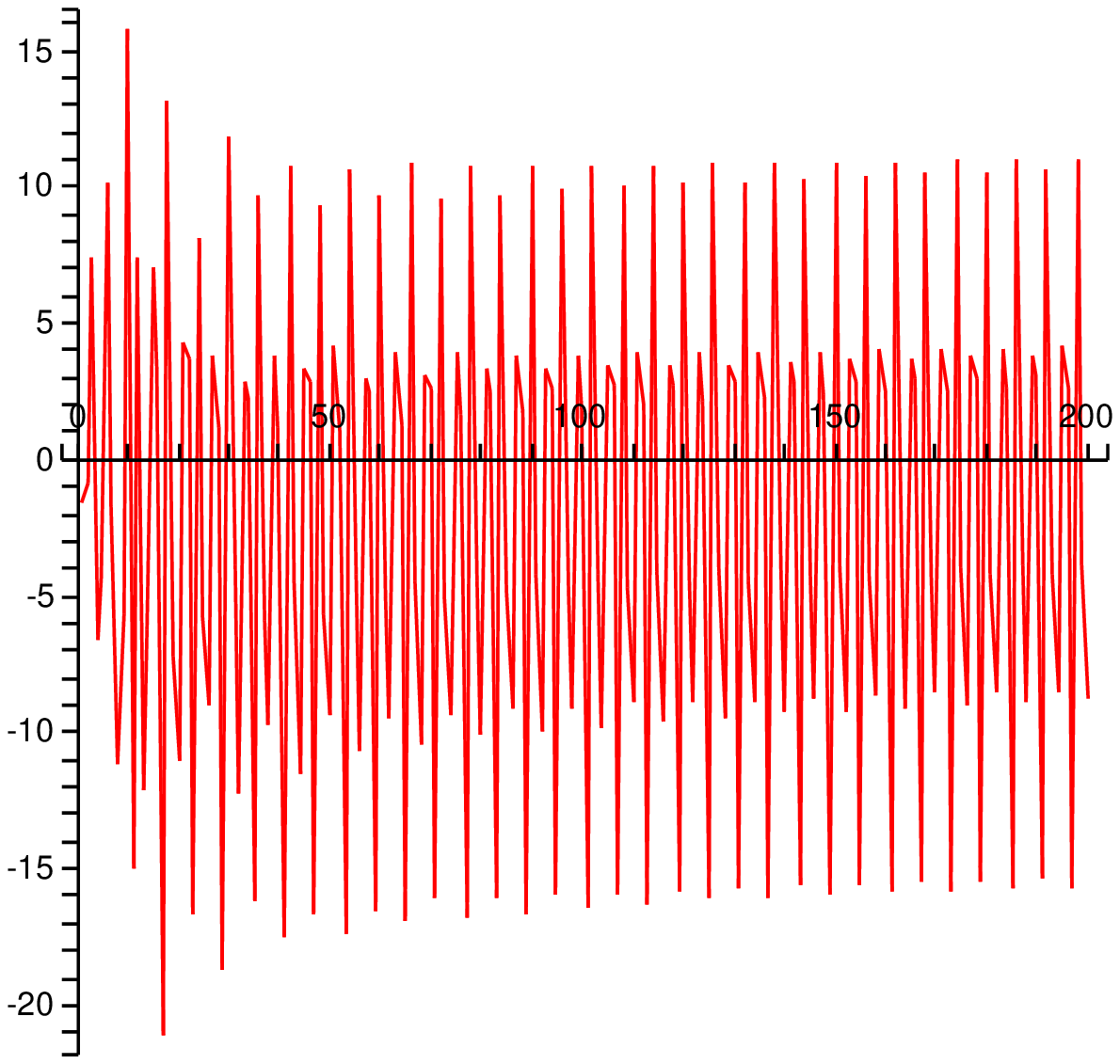}
\hspace*{0.5 truecm}
\includegraphics[width=50truemm]{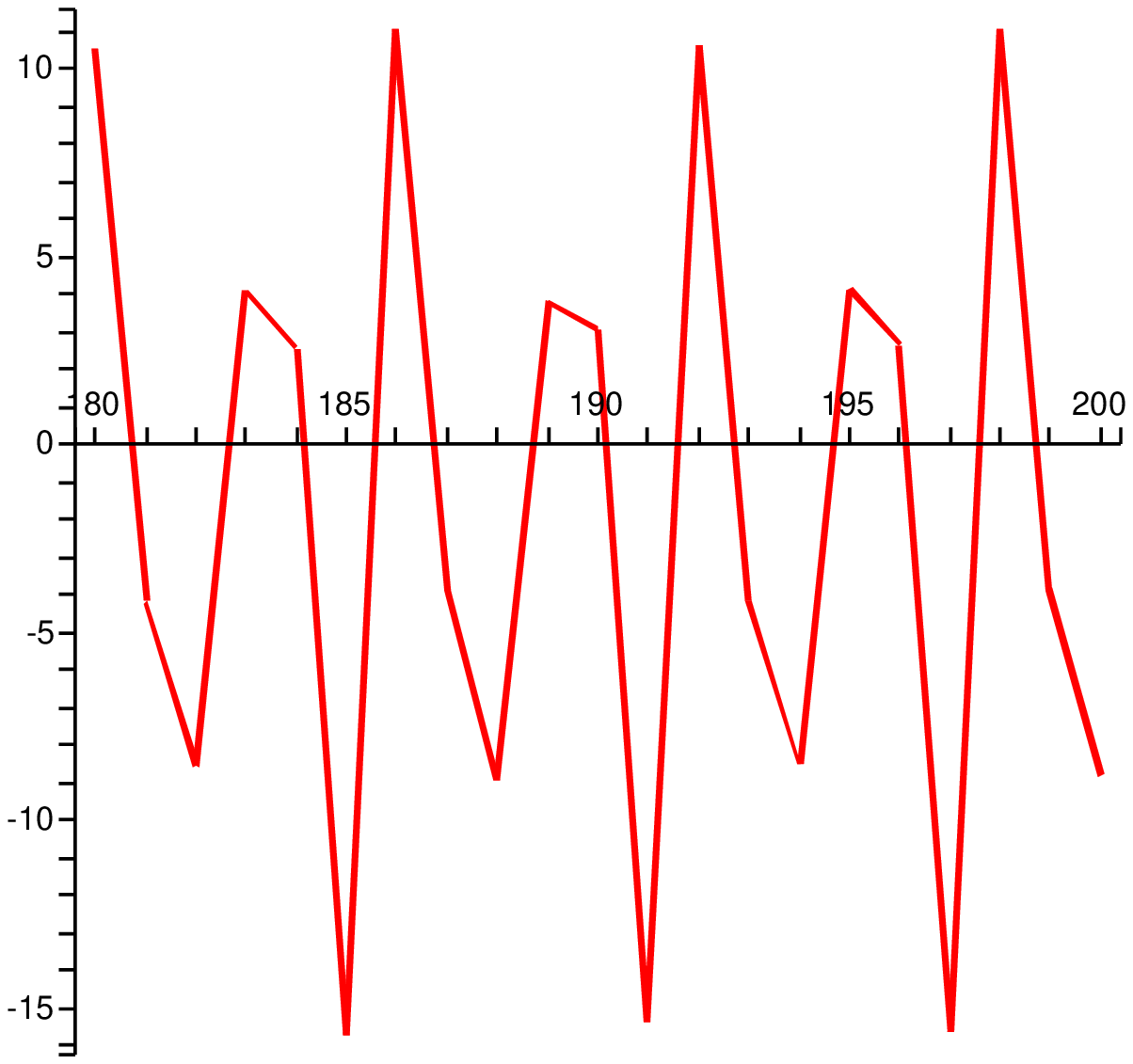}}
\end{center}

\vspace*{-0.8 truecm}

\caption{\label{error-fig}
Permutations with distinct cycle lengths: the 
approximation error as measured by $n^3(f_n/\hat f_n^{[1]}-1)$,
with $\hat f_n^{[1]}$ truncated after $n^{-3}$ terms,
for $n=1\,.\,.\,200$ (left) and for $n=180\,.\,.\,200$ (right).}
\end{figure}

The terms given in~\eqref{at1} provide quite a good approximation. 
Figure~\ref{error-fig} displays a comparison between $f_n$ and its asymptotic
approximation $\hat f_n^{[1]}$, up to terms of order $n^{-3}$.
We find
\[
f_n = \hat f_n^{[1]}\left(1+\frac{R_n}{n^3}\right),\qquad
|R_n|\le 22,
\]
for all $1\le n\le 1000$, and expect the bound to remain valid for all
$n\ge 1$.

\smallskip

\emph{Expansions at $z=-1$ and at $z=\omega,\omega^2$.}
We shall content ourselves with brief indications on the shape of 
the corresponding singular expansions. 
Note that Figure~\ref{error-fig} clearly indicates  the presence of a term of the form
\[
\frac{\Omega(n)}{n^3}
\]
in the full asymptotic expansion of $f_n$, where $\Omega(n)$ is a function with period 6. 
This motivates
an investigation of the behaviour of $f(z)$ near primitive square and cube roots of unity.

Start with $z\to-1^{+}$. The definition of $f(z)$ implies that
\[
f(z)\mathop{\sim}_{z\to-1^+} (1+z)\prod_{k=2}^\infty \left(1+\frac{(-1)^k}{k}\right)=(1+z)
\]
(the  infinite   product telescopes).    Set $\zz=1+z$   and  restrict
attention to the type of the expansions at $z=-1$. Only half of the
polylogarithms in~\eqref{chain2} are singular, so that the expansion at $z=-1$ is 
of the \emph{type}
\[
\zz \exp\left(\zz\log\zz +\zz +\zz^2+\zz^3\log\zz+\cdots\right)
=
\zz+\zz^2\log\zz +\zz^2+\zz^3\log^2\zz+\cdots\,.
\]
There, we have replaced all unspecified coefficients by the constant~1 for readability.
This singular form results in a contribution to the asymptotic form of $f_n$:
\[
\hat f_n^{[-1]}=(-1)^n\left(\frac{d_3}{n^3}+\frac{d_{4,1}\log n+d_{4,0}}{n^4}+\cdots\right).
\]
(Compared to roots of unity of higher order,
the case $z=-1$  is special, because  of the factor $(1+z)$
explicitly present in the definition of $f(z)$.)
A simple calculation shows that $d_3=2$, so that
\[
\hat f_n^{[-1]}\sim 2\frac{(-1)^n}{n^3}.
\]

Next, let $\omega=e^{2i\pi/3}$ and set $\zz:=(1-z/\omega)$. The 
type of the expansion at $z=\omega$
is
\[
\begin{array}{lll}
f(z)&\sim &\ds f(\omega)\exp\left(\zz+\zz^2\log\zz+\zz^3+\cdots\right)
\\
&\sim&\ds f(\omega)\left(1+\zz+\zz^2\log\zz+\zz^2+\zz^3\log\zz+
\zz^4\log^2\zz+
\cdots\right),
\end{array}
\]
since now every third polylogarithm is singular at $z=\omega$.
This induces a contribution of the form
\[
\hat f_n^{[\omega]}=\omega^{-n}\left(\frac{e_3}{n^3}+\frac{e_{4,1}\log n+e_{4,0}}{n^4}+\cdots\right),
\]
arising from $z=\omega$,  and similarly, for  a conjugate contribution
arising from $\omega^2$.  Another simple calculation shows that
\[
e_3=3f(\omega),\]
and leaves us with the task of estimating $f(\omega)$. The use of the formula,
\[
\prod_{k\ge1}\left(1+\frac{s}{n+a}\right)e^{-s/n}
=\frac{\Gamma(1+a)e^{-s\gamma}}{\Gamma(1+s+a)},
\]
a mere avatar of the product formula for the Gamma function, yields then
easily
\[
f(\omega)=
\frac{3\Gamma(\frac23)}{\Gamma(\frac13+\frac{\omega}{3})
\Gamma(\frac23+\frac{\omega^2}{3})}.
\]

The fluctuations of period~6 evidenced by Figure~\ref{error-fig} are 
thus fully explained: one has
\[
\hat f_n^{[-1]}+\hat f_n^{[\omega]}+\hat f_n^{[\omega^2]}=\frac{\Omega(n)}{n^3}+
O\left(\frac{\log n}{n^4}
\right),
\]
where the periodic function $\Omega$ is $(\Re$ designates a real part)
\[
\Omega(n)=2(-1)^n +2\Re\left(\omega^{-n}
\frac{9\Gamma(\frac23)}{\Gamma(\frac16+\frac{i\sqrt{3}}{6})
\Gamma(\frac12-\frac{i\sqrt{3}}{6})}\right).
\]

\smallskip

\emph{Expansions at $z=\eta$, a primitive $\ell$th root of unity.}
Let  $\eta=\exp(2i\pi/\ell)$ and  $\zz=(1-z/\eta)$.   The expansion of
$f(z)$ is now of the type
\[
f(z)\sim f(\eta)\exp\left(\zz+\cdots+\zz^{\ell-1}\log\zz+\zz^{\ell-1}+\cdots\right),
\]
where~$\zz^{\ell-1}\log\zz$ corresponds to the singular term in $\Li_{\ell}(z^\ell)$.
Consequently,  fluctuations start appearing at  the  level of terms of
order $n^{-\ell}$ in the asymptotic expansion of~$f_n$ as $n\to+\infty$.
The value of $f(\eta)$ is expressible in terms of Gamma values at algebraic points,
as we have seen when determining $f(\omega)$. The coefficients 
in the expansion also involve values of the $\psi$-function
($\psi(s)=\frac{d}{ds}\log\Gamma(s)$) and its derivatives
at rational points, which include $\zeta$ values as particular cases.

\begin{figure}

\begin{center}
\hbox{\hspace*{1.5truecm}\includegraphics[width=50truemm]{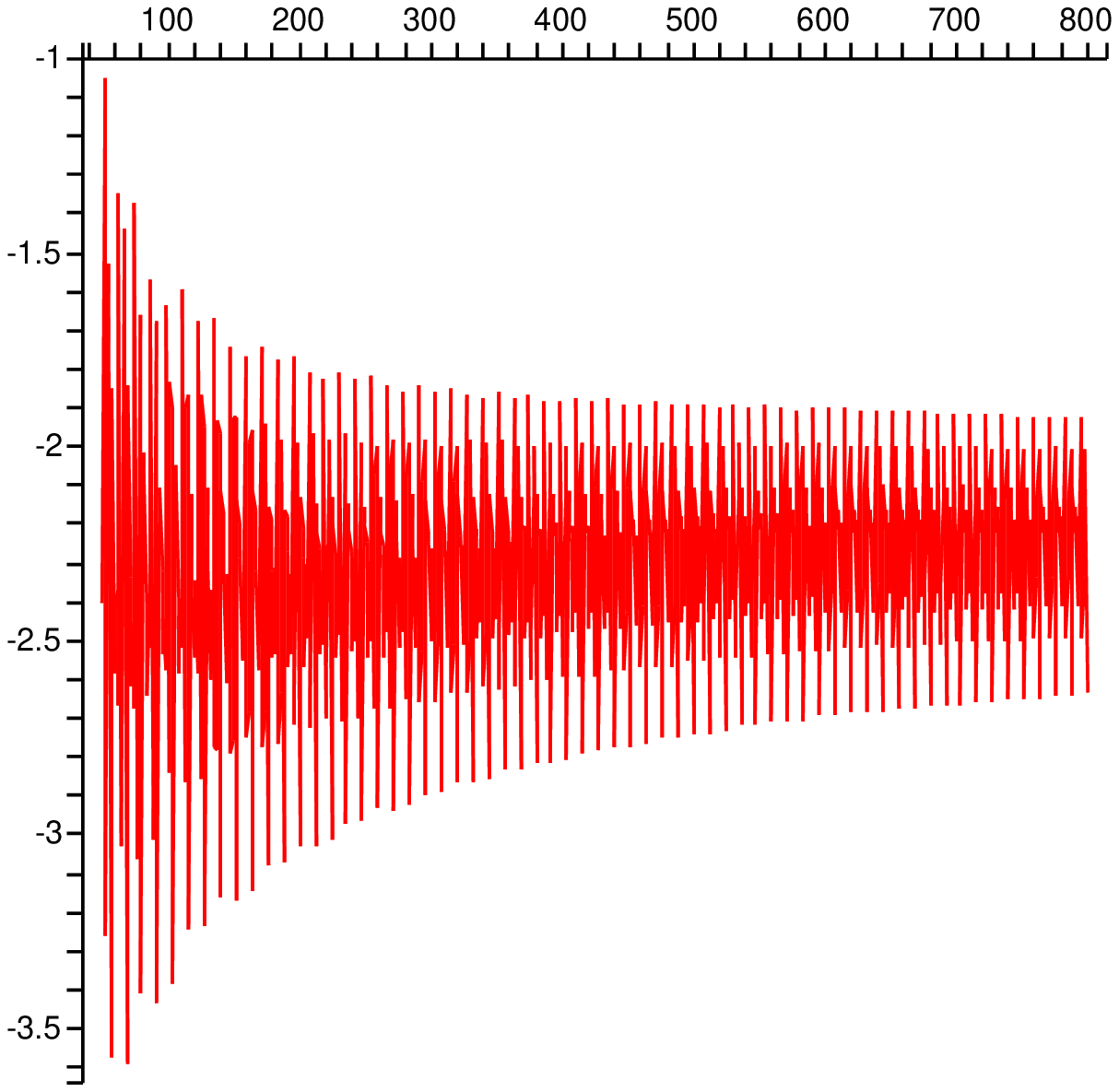}
\hspace*{0.5truecm}\includegraphics[width=50truemm]{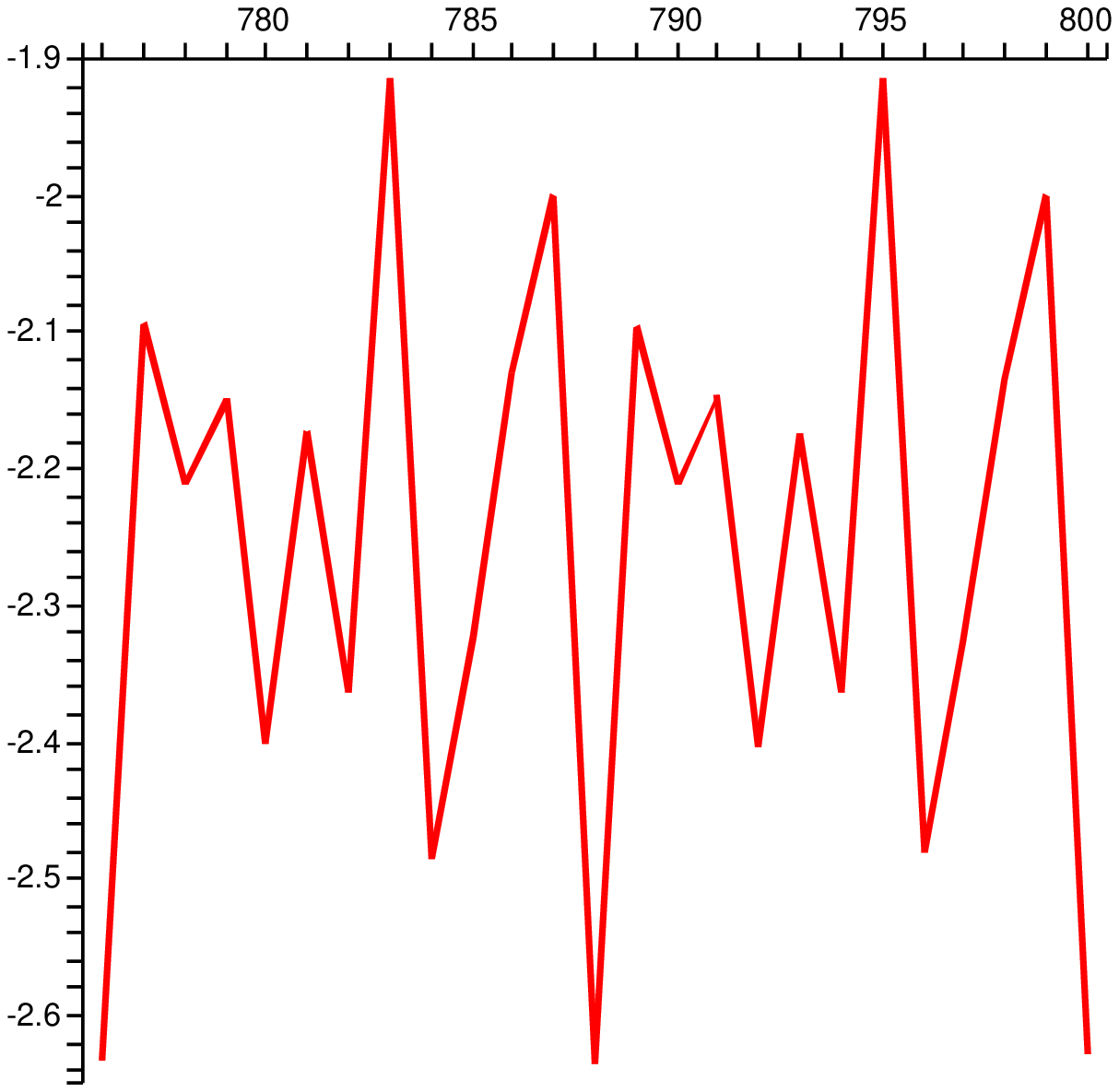}}
\end{center}

\vspace*{-0.8 truecm}

\caption{\label{error4-fig}
Permutations with distinct cycle lengths: the approximation error as measured
by
$n^4\log^3n (f_n/(\hat f_n^{[1]}+\hat f_n^{[-1]}
+\hat f_n^{[\omega]}+\hat f_n^{[\omega^2]})-1)$,
with the $\hat f_n$ truncated after $n^{-3}$ terms,
for $n=50\,.\,.\,800$ (left) and for $n=776\,.\,.\,800$ (right).}

\end{figure}

\begin{proposition} \label{grkn-prop}
The  probability that  a permutation is made  of cycles of  distinct
lengths admits a full asymptotic expansion of the form
\[
\begin{array}{lll}
f_n &\sim& \ds
 e^{-\gamma}+\frac{e^{-\gamma}}{n}+\frac{e^{-\gamma}}{n^2}\left(-\log n-1-\gamma+\log 2\right)\\
&&\ds 
+\frac{1}{n^3}
\left[e^{-\gamma}\left(\log^2n+c_{3,1}\log n+c_{3,0}\right)
+2(-1)^n+2\Re\left(
\hbox{$\frac{9\Gamma(\frac23)\omega^{-n}}{\Gamma(\frac16+\frac{i\sqrt{3}}{6})
\Gamma(\frac12-\frac{i\sqrt{3}}{6})}$}\right)\right]
\\
&&\ds {}+ \sum_{r\ge 4}\frac{P_r(n)}{n^r},
\end{array}
\]
with $c_{3,1}$ and $c_{3,0}$ as given by~\eqref{c30eq}.
There, $P_r(n)$ is a polynomial of exact degree~$r-1$ in $\log n$ with
coefficients that  are  periodic functions  of  $n$  
with period $D(r)=\lcm(2,3,\ldots,r)$.
\end{proposition}

Figure~\ref{error4-fig}  displays   the error  of   the  approximation
obtained  by incorporating all terms  till  order $n^{-3}$ included in
the asymptotic expansion of  $f_n$.  Fluctuations of period~12 (due  to
the additional presence  of $i=\sqrt{-1}$) start  making  an appearance.

\section{Permutations, polynomials, and trees}\label{hodgepodge-sec}

We now examine several combinatorial problems 
related to permutations, polynomials over finite fields, and trees that 
are amenable to the hybrid method. The detailed treatment of 
permutations with distinct cycle lengths can serve as a beacon
for the analysis of similar infinite product generating functions,
and accordingly our presentation of each example will be quite succinct.

In  the examples that  follow, the function~$f$ whose coefficients are
to be analysed is such that there is an increasing  family of sets $Z^{(1)},Z^{(2)},\ldots$
(ordered  by inclusion, and with elements being roots of
unity),  attached to   a collection  of  asymptotic expansions  having
smaller and smaller error  terms. In  that  case,  a full asymptotic
expansion  is  available  for  the   coefficients of~$f$.  The general
asymptotic shape  of  $[z^n]f$  involves standard   terms  of the  form
$n^{-p}\log   ^q  n$   modulated   by   complex   exponentials,  since
$[z^n]\sigma(z/\zeta)=\zeta^{-n}\sigma_n$. We formalize this notion by a definition.

\begin{definition}
A sequence $(f_n)$ is said to admit a \emph{full asymptotic expansion with 
oscillating coefficients} if it is of the form
\[
f_n \sim \sum_{r\ge0} \frac{P_r(n)}{n^{\alpha_r}}
\qquad(n\to\infty),
\]
where the exponents~$\alpha_r$ increase to $+\infty$ and 
each $P_r(n)$ is a polynomial 
in~$\log n$ whose coefficients are periodic functions of~$n$.
\end{definition}

In a small way, Proposition~\ref{grkn-prop} and  the forthcoming statements,
Propositions~\ref{pi2-prop}--\ref{for-prop},
can  be regarded  as  analogues, in
the  realm of functions of  slow growth near  the unit  circle, of the
Hardy-Ramanujan-Rademacher
 analysis~\cite{Andrews76,Ayoub63,Hardy78} of partition generating functions,
the latter exhibiting a very fast growth (being of infinite order) as $|z|\to1$.

\subsection{Permutations admitting an $m$-th root}\label{mthroot-sec}

%
%
%
%
%

The problem of determining the number of permutations that are squares
or  equivalently  ``have a  square   root'' is   a  classical  one  of
combinatorial       analysis:    see       Wilf's   vivid      account
in~\cite[\S4.8]{Wilf94}. The problem admits an obvious generalization.
We  shall let  $\Pi_m(z)$ be  the  exponential generating  function of
permutations that are $m$th powers  or, if one  prefers, admit an  $m$th
root.

\parag{How many permutations have square roots?}
For the generating function $\Pi_2$, we follow Wilf's account. 
Upon squaring a permutation~$\tau$, each cycle of even length of~$\tau$
falls apart into two cycles of half the length, while an odd cycle 
gives rise to a cycle of the same length. Hence, if~$\sigma$ has a square root, then the number
of cycles it has of each even length must be even. By general principles of
combinatorial analysis, the exponential generating function 
of the number $\Pi_{2,n}$ of permutations of $n$ elements that have square roots satisfies
\begin{equation}\label{bender0}
\begin{array}{lll}
\ds \Pi_{2}(z)\,:=\, \sum_{n\ge0} \Pi_{2,n}\frac{z^n}{n!} &=&
\ds e^{z}\cosh(z^2/2)e^{z^3/3}\cosh(z^4/4)e^{z^5/5}\cdots \\
&=&\ds \exp\left(\frac12\log\frac{1+z}{1-z}\right)\prod_{m\ge1}\cosh \left(\frac{z^{2m}}{2m}\right)\\
&=& \ds \sqrt{\frac{1+z}{1-z}}\prod_{m\ge1}\cosh \left(\frac{z^{2m}}{2m}\right).
\end{array}
\end{equation}
The series starts as
\[
\Pi_2(z)=1+z+\frac{z^2}{2!}+3\frac{z^3}{3!}+12\frac{z^4}{4!}+60\frac{z^5}{5!}+
270\frac{z^6}{6!}+\cdots\,,\] its  coefficients   being~\EIS{A003483},
with  the quantity    $[z^n]\Pi_2(z)$   representing the probability    that  a
permutation  is  a  square. The   GF  of~\eqref{bender0} is   given by
Bender~\cite[p.~510]{Bender74}         who     attributes     it     to
Blum~\cite{Blum74}.  It is interesting to note that Bender mentions
the following estimate from~\cite{Blum74}, 
\begin{equation}\label{pi21}
[z^n]\Pi_{2}(z) \sim \frac{2}{\sqrt{\pi n}} e^G,
\end{equation}
and derives it  by   an  application of   the Tauberian  theorem\footnote{%
	Regarding Tauberian side conditions, B\'ona, McLennan, and White~\cite{BoMcWh00}
	prove by elementary combinatorial arguments
	 that the sequence~$\Pi_{2,n}$ is monotonically nonincreasing
	in~$n$.
}
of  Hardy,
Littlewood, and Karamata.  Accordingly,  no error terms  are available,
given the nonconstructive character of classical Tauberian theory.  We
state:

\begin{proposition} \label{pi2-prop}
The probability that a random permutation of size~$n$ has a square-root
admits a full asymptotic expansion with oscillating coefficients.
In particular, it satisfies
\begin{equation}\label{pif}
[z^n]\Pi_{2}(z) \sim \sqrt{\frac{2}{\pi n}} e^G
\left[1-\frac{\log n}{n}+\frac{c_3+(-1)^n}{4n}\right]
-\frac{2e^G(-1)^{\lfloor n/2\rfloor}}{n^2}
+O\left(\frac{\log n}{n^{5/2}}\right),
\end{equation}
where
\begin{equation}\label{pi22}
\left\{
\begin{array}{lll}
e^G&=&\ds \prod_{k\ge1}\cosh \left(\frac{1}{2k}\right)\quad \doteq\quad  1.22177\,95151\,92536
\\
c_2&=& \ds \sum_{k\ge1} \left(\frac{1}{2k}-\tanh\big(\frac{1}{2k}\big)\right), \qquad
c_3=-12+16\log 2+4\gamma+2c_2.
\end{array}
\right.
\end{equation}
\end{proposition}
\begin{proof} First a rough analysis suffices to see that
\[
\Pi_2(z) \mathop{\sim}_{z\to1} \sqrt{\frac{2}{1-z}}e^G \qquad\hbox{and}\qquad
\Pi_2(z) \mathop{\sim}_{z\to-1} \sqrt{\frac{1+z}{2}}e^G .
\]
In order to refine these expansions, introduce the normalized tangent numbers
by
\[
\tan z=\sum _{m\ge 0}\tau_m z^{2m+1},
\qquad\hbox{so that}\quad
\log\cos(z)=\sum_{m\ge1}\tau_{m-1} \frac{z^{2m}}{2m}.
\]
The usual exp-log reorganization of the $\Pi_2$ series yields
\[
\Pi_2(z)=\sqrt{\frac{1+z}{1-z}}\exp\left(\sum_{m\ge1} \frac{(-1)^{m-1}}{m2^{2m+1}}
\tau_{m-1}\Li_{2m}\left(z^{4m}\right)\right).\]
In passing, this provides for~$G$, the fast convergent series
\[
G=\sum_{m\ge1} \frac{(-1)^{m-1}}{m2^{2m+1}}\tau_{m-1}\zeta(2m),
\]
on which the numerical estimate of~\eqref{pi22} is based.

Next, take out the $e^G$ factor, leading to
\begin{equation}\label{pi23}
\Pi_2(z)=e^G \sqrt{\frac{1+z}{1-z}}
\exp\left(\sum_{m\ge1} \frac{(-1)^{m-1}}{m2^{2m+1}}\tau_{m-1}
\left[\Li_{2m}\left(z^{4m}\right)-\zeta(2m)\right]\right).
\end{equation}
At~$z=1$, the largest singular term in the exponential arises from~$\Li_2(z^4)$, the 
contributions from the other polylogarithms being of smaller order 
\[
\Pi_2(z) \mathop{\sim}_{z\to1}
 e^G\sqrt{\frac{2}{1-z}}\left(1+2(1-z)\log(1-z)+c^\star (1-z)+
O\left((1-z)^2\log^2(1-z)\right)\right),
\]
where, here and later, $c^\star$ designates a computable constant
that we leave unspecified for the purpose of readability.
Similarly, at $z=-1$, we find 
\[ 
\Pi_2(z) \mathop{\sim}_{z\to-1} e^G\sqrt{\frac{1+z}{2}}
\left(1+O\left((1+z)\log(1+z)\right)\right).
\]
At $z=i=\sqrt{-1}$ (hence at $z=-i$, by conjugacy), 
we have

\[
\Pi_2(z)\mathop{\sim}_{z\to i} \frac{1+i}{2}\left(P^\star(1-z/i)+2(1-z/i)\log(1-z/i)
+O((1-z/i)^{3/2})\right),
\]
with $P^\star$ an unspecified polynomial that does not leave a trace in the coefficients' expansion. 

The singular contribution to $\Pi_2$ 
arising from any root of unity~$\zeta$ of order $\ge8$ 
is at most  $(O(1-z/\zeta)^3\log(1-z/\zeta))$, which translates to an
$O(n^{-4})$ term. The proof of~\eqref{pif} is then  completed
upon making $c^\star$ in the expansion at~1 explicit,
which introduces the new constants $c_2,c_3$.
Existence of the full expansions finally follows from the usual analysis of polylogarithms 
of powers, taken at roots of unity.
\end{proof}

\parag{How many permutations have an $m$th root?}
Like the previous one, this problem is briefly mentioned in Bender's survey~\cite{Bender74}.
We follow again Wilf's exposition~\cite[\S4.8]{Wilf94}.
For a pair $\ell,m$ of positive integers, we 
define $\ggcd{\ell,m}$ to be
\[
\ggcd{\ell,m}:=\lim_{j\to\infty} \gcd(\ell^j,m).
\]
(Thus, $\ggcd{\ell,m}$ gathers from the prime decomposition
of~$m$  all  the factors  that  involve  a prime   divisor of~$\ell$.) The
characterization   of    permutations  that are  $m$th    powers  then
generalizes~\cite{Bender74,Pouyanne02,Turan70,Wilf94}:
a permutation has an $m$th root if and only if, for each~$\ell=1,2,\ldots$, 
it is true that 
the number of cycles of length~$\ell$ is a multiple of~$\ggcd{\ell,m}$.
This observation leads to an expression for the
corresponding generating function. Indeed, define the ``sectioned exponential'',
\[
\exp_d(z):=\sum_{n\ge0} \frac{x^{dn}}{(dn)!}
=\frac{1}{d}\sum_{j=0}^{d-1}\exp\left(e^{2ij\pi/d}z\right),
\]
so that $\exp_1(z)=e^z$ and $\exp_2(z)=\cosh(z)$. The exponential generating function of 
permutations that are $m$th powers is then
\begin{equation}\label{pim}
\Pi_{m}(z)=\prod_{\ell=1}^\infty \exp_{\ggcd{\ell,m}}\left(\frac{z^\ell}{\ell}\right).
\end{equation}

The generating function of~\eqref{pim} has been investigated by Pouyanne~\cite{Pouyanne02},
whose paper provides the first order asymptotic estimate of $[z^n]\Pi_m$.
There is a fundamental factorization,
\begin{equation}\label{pim1}
\Pi_m(z)=A_m(z)\cdot B_m(z),
\end{equation}
where $A_m$ gathers from the product~\eqref{pim} the numbers $\ell$ 
that are relatively prime to~$m$ and $B(z)$ gathers the rest.

The factor $A_m$ is found by series rearrangements to be
an algebraic function expressible by radicals,
\[
A_m(z)=\prod_{k~|~m} (1-z^k)^{-\mu(k)/k},\] with $\mu(k)$ the   M\"obius
function.    For $m=2$    (square   permutations),   this   is     the
ubiquitous prefactor~$A_2(z)=\sqrt{(1+z)/(1-z)}$. For $m=6$, for instance, 
the prefactor becomes
\[
A_6(z)=\left(\frac{1+z}{1-z}\right)^{1/3}\left(\frac{1+z+z^2}{1-z+z^2}\right)^{1/6}.\]

The factor $B_m$ is a transcendental function that admits the unit circle 
as a natural boundary~\cite{Pouyanne02}. It is expressible as an infinite product
of sectioned exponentials:
\begin{equation}\label{pim3}
B_m(z)=
\prod_{{\ell=1}\atop
{\gcd(\ell,m)>1}}^\infty \exp_{\ggcd{\ell,m}}\left(\frac{z^\ell}{\ell}\right).
\end{equation}
For $m=2$, this is the infinite product of hyperbolic cosines. For $m=6$, one has
\[
B_6(z)=
\prod_{\ell\equiv 0 \mod 6} \exp_6\left(\frac{z^\ell}{\ell}\right)
\prod_{\ell\equiv 2, 4\mod 6} \exp_2\left(\frac{z^\ell}{\ell}\right)
\prod_{\ell\equiv 3 \mod 6} \exp_3\left(\frac{z^\ell}{\ell}\right).
\]
These singular factors can be analysed just like in the case of~$\Pi_2$
by an exp-log transformation. 
One first observes that the limit value $B_m(1)$ is well defined,
since the infinite product converges at least as fast as $\prod(1+\ell^{-2})$.
It is seen next that singularities 
are at roots of unity, and the radial expansions can be computed 
in the usual way from the polylogarithmic expansion. We can now state
a (somewhat minor) improvement over~\cite{Pouyanne02}:

\begin{proposition} The probability that a random permutation of size~$n$ has
an~$m$th root admits a full asymptotic expansion with oscillating
coefficients. To first asymptotic order, it satisfies
\[
[z^n]\Pi_m(z)\sim \frac{\varpi_m}{n^{1-\varphi(m)/m}},
\qquad
\varpi_m:= \frac{B_m(1)}{\Gamma(\varphi(m)/m)}\prod_{k~|~m} k^{-\mu(k)/k},
\]
where $\varphi(m)$ is the Euler totient function and
\begin{equation}\label{pouyanne}
B_m(1)=
\prod_{{\ell=1}\atop
{\gcd(\ell,m)>1}}^\infty \exp_{\ggcd{\ell,m}}\left(\frac{1}{\ell}\right).
\end{equation}
In particular, when $m=p$ is a prime number, one has
\begin{equation}\label{turan}
[z^n]\Pi_p(z)\sim \frac{\varpi_p}{n^{1/p}},
\qquad
\varpi_p:=\frac{p^{1/p}}{\Gamma(1-1/p)}\prod_{\ell=1}^\infty \exp_p\left(\frac{1}{\ell p}\right).
\end{equation}
\end{proposition}
For instance, we find:
\[
[z^n]\Pi_6(z)=B_6(1)\left[ \frac{\root6\of{12}}{\Gamma(\frac13)}n^{-2/3}
+\frac{2\root12\of{12}}{\Gamma(\frac16)}
\Re\left({e^{-i\frac{\pi}{3}(n+\frac14)}}\right)n^{-5/6}+
O\left(\frac1n\right)\right],\] 
and~\eqref{turan} improves  upon an
early estimate of Tur\'an~\cite[Th.~IV]{Turan70}.

\subsection{Pairs of permutations having the same  cycle type}
Given a permutation $\sigma$, its \emph{cycle type}
is the  (unordered) multiset formed with the 
lengths of the cycles entering its decomposition into cycles.
For a permutation of size~$n$, this type can be equivalently represented 
by a partition of the integer~$n$: 
for instance $(2^3,5,7^2)$ represents the profile of
any permutation of size~25 that has three cycles of length~2, one cycle of length~5, and
two cycles of length~7. The probability $a_{\pi}$ that a random permutation
of size $n$ has profile 
\[
\pi=(1^{n_1},2^{n_2},\ldots), \qquad \hbox{with}\quad
1n_1+2n_2+\cdots =n
\]
is, by virtue of a well-known formula~\cite[p.~233]{Comtet74},
\[
a_{\pi}= \prod_{i\ge 1} \frac{1}{i^{n_i} n_i!},
\]
corresponding to the generating function
in infinitely many variables
\begin{equation}\label{cauchyy}
\Phi(z;x_1,x_2,\ldots)=\exp\left(x_1z+x_2\frac{z^2}{2}+x_3\frac{z^3}{3}+\cdots \right),
\end{equation}
which is such that $a_\pi=[z^nx_1^{n_1}x_2^{n_2}\cdots]\Phi$.

In this subsection, we estimate \emph{the probability that two permutations of 
size~$n$ taken uniformly and independently at random have the same cycle type.}
(Each pair of permutations is taken with probability $1/n!^2$.)
The quantity to be estimated is thus
\[
W_n = \sum _{\pi\, \vdash\,  n} \big(a_\pi\big)^2,
\]
where      the          summation        ranges       over         all
partitions~$\pi=(1^{n_1},2^{n_2},\ldots)$ of~$n$.    This  problem was
suggested to   the   authors by  the   reading  of  a short   note  of
Wilf~\cite{Wilf06},   who     estimated  the   probability   that  two
permutations have the same number of cycles  (the answer to the latter
question turns out to be asymptotic to $1/(2\sqrt{\pi\log n})$).

Given~\eqref{cauchyy}, it is not hard to find the generating function of the sequence~$(W_n)$:
\begin{equation}\label{wi1}
W(z):=\sum_{n\ge0}W_n z^n=\prod_{k=1}^\infty I\left(\frac{z^k}{k^2}\right)
\qquad\hbox{where}\quad
I(z)=\sum_{n\ge0} \frac{z^n}{n!^2},
\end{equation}
the reason being that
\[
W(z)=\sum_{n\ge0}\sum_{\pi\,\vdash\,  n} \frac{z^{n_1+2n_2+\cdots}}{\prod_i i^{2n_i} n_i!^2}
=\prod_{i\geq 1}\left(\sum_{n_i\geq0} \frac{z^{i n_i}}{i^{2n_i} n_i!^2}\right).
\]
The function written~$I(z)$, which is obviously entire,
is a variant of the Bessel function~$I_0$ 
(see, e.g., ~\cite{WhWa27}): $I(z)=I_0(2\sqrt{z})$.
Also, from~\eqref{wi1}, the expansion of $W(z)$ is readily computed: one has
\[
W(z)=1+z+2\frac{z^2}{2!^2}+14\frac{z^3}{3!^2}+146\frac{z^4}{4!^2}+
2602\frac{z^5}{5!^2}+\cdots ,
\]
where the coefficients are \EIS{A087132} (``sum of the squares of the sizes of
conjugacy classes in the symmetric group$~\mathfrak{S}_n$'').

\begin{proposition}
The probability~$W_n$ that two permutations of
size $n$ have the same cycle type satisfies
\begin{equation}
W_n= \frac{W(1)}{n^2}+O\left(\frac{\log n}{n^3}\right),
\qquad W(1)=\prod_{k\geq 1}I\left(\frac{1}{k^2}\right)\doteq 4.26340\,35141\,52669.
\end{equation}
Furthermore, this probability admits a full asymptotic expansion with oscillating coefficients.
\end{proposition}
\begin{proof}
We shall only sketch the analysis of the dominant asymptotic term, the rest being 
by now routine.
From the expression of $W(z)$,
the exp-log transformation yields
$$
W(z)=\exp\left(\sum_{k\geq1}H\left(\frac{z^k}{k^2}\right)\right)
\qquad \mathrm{where}\quad
H(z):=\log\left(I(z)\right). 
$$  
Let $h_\ell=[z^\ell]H(z)$,
that is, $H(z)=\sum_\ell h_\ell z^\ell=z-\frac{1}{4}z^2+\cdots$
(the sequence $(h_\ell\ell!^2)$ is \EIS{A002190}, which occurs in the enumeration of 
certain pairs of permutations by Carlitz).
Then 
$$
W(z)=\exp\left(\sum_{k\geq 1}\sum_{\ell\geq 1}h_\ell\frac{z^{k\ell}}{k^{2\ell}}\right)
=\exp\left( \sum_{\ell\geq 1}h_\ell \mathrm{Li}_{2\ell}(z^\ell)\right).
$$
This expression ensures that the hybrid method can be applied at any
order, with $W(z)$ being of global order $a=0$. The
first term of the asymptotic estimate is provided by the factorization
relative to $Z=\{1,-1\}$, namely $W(z)=P(z)Q(z)$, with
\[
P(z)=\exp\left(\mathrm{Li}_2(z)-\frac{1}{4}\mathrm{Li}_4(z^2)\right),\qquad 
Q(z)=\exp\left(\sum_{\ell\geq 3}h_\ell\mathrm{Li}_{2\ell}(z^\ell)\right).
\]
Since $Q(z)$ is clearly $\mathcal{C}^4$ on the closed unit disk, the
hybrid method applies: with the notations of
Theorem~\ref{hybrid2-thm}, one can take $s=4$, to the effect that 
 $u_0=2$. Using the algorithmic scheme of
Section~\ref{discyc-sec}, we find 
\[\renewcommand{\arraycolsep}{2truept}
\begin{array}{lll}
W(z) &\ds \mathop{=}_{z\to
1^-} & \ds W(1)\left[1-(1-z)L+\frac{1}{2}(1-z)^2L+\frac{1}{2}(1-z)^2L^2+O((1-z)^2)\right]
\\
W(z) & \ds \mathop{=}_{z\to -1^+} & \ds W(-1)\left[1+O((1+z)^2)\right],
\end{array}
\]
with $L\equiv L(z)=\log(1-z)^{-1}$.
Theorem~\ref{hybrid2-thm} then directly yields $W_n\sim W(1)/n^2$, and
further asymptotic terms can easily be extracted. 
\end{proof}

\smallskip

Regarding other statistics on  \emph{pairs} of permutations, it is well
worth mentioning Dixon's recent  study   of the probability that   two
randomly chosen permutations generate a transitive group~\cite{Dixon05}.
For the symmetric group~$\mathfrak{S}_n$, this probability is 
found to be asymptotic to 
\[
1-\frac{1}{n}-\frac{1}{n^2}-\frac{4}{n^3}-\frac{23}{n^4}-\cdots\,\]
and, up to exponentially smaller order terms, this
expansion also gives the probability
that two random permutations generate the whole symmetric group.
In that case, the analytic engine is Bender's theory of
coefficient extraction in divergent series~\cite{Bender74}.

\subsection{Factorizations of polynomials over finite fields}

Factoring polynomials in~$\Q[X]$ is a problem of interest
in symbolic computation, as it has implications in the determination
of partial fraction expansions and symbolic integration, for instance.
Most of the existing algorithms proceed by a reduction to
a few factorizations of polynomials with coefficients in some finite field~$\F_q$
(with~$q$ a prime power): see the books by Berlekamp, Knuth, or von zur Gathen and Gerhard~\cite{Berlekamp68,GaGe99,Knuth98}.
Several of the algorithms for factoring polynomials over finite fields
involve what is known as the \emph{distinct degree factorization}. To wit,
let $\phi\in \F_q[X]$ be completely factored as
\[
\phi = \prod_{j} \iota_j ^ {n_j},\]
where the $\iota_j$ are distinct irreducible polynomials. The distinct degree
factorization produces the decomposition of the square-free part of $\phi$ under
the form
\[
\operatorname{squarefree}(\phi)=\beta_1 \beta_2 \cdots \beta_n,
\qquad\hbox{where}\quad 
\beta_k = \prod_{\deg(\iota_j)=k} \iota_j\,.
\]
We refer to the paper~\cite{FlGoPa01} by three of us for details.

The distinct degree factorization of~$\phi$
as defined above coincides with the complete factorization of~$\phi$
if and only $\phi$ has all its irreducible factors of 
different degrees. 
For the purpose of analysis of algorithms,
it is then of interest to quantify the probability of that event.
This problem has been considered independently by the authors of~\cite{FlGoPa96,FlGoPa01}
and by Knopfmacher \& Warlimont in~\cite{KnWa95}.

First, without loss of generality, we restrict attention to monic polynomials. If~$q$ is the cardinality of the base field $\Fq$, the number of monic polynomials 
of degree $n$ is $P_n=q^n$, with ordinary generating function
\[
P(z)=\sum_{n\ge0} q^n z^n = \frac{1}{1-qz}.\]
The unique factorization property implies that the class 
$\cal P$ of all (monic) polynomials is isomorphic to the 
class obtained by applying the Multiset construction to
the class $\cal I$ of (monic) irreducible polynomials. 
This, by well known principles of combinatorial analysis~\cite{Berlekamp68,FlSe06}
translates to an equation binding the generating function $I(z)$ of irreducible polynomials,
namely,
\begin{equation}\label{gau}
P(z)=
\prod_{n=1}^\infty (1-z^n)^{-I_n}
=\exp\left(I(z)+\frac12I(z^2)+\frac{1}{3}I(z^3)+\cdots\right),
\end{equation}
where~$I_n$ is the number of monic irreducible polynomials of degree~$n$.
The solution is obtained by taking a logarithm 
of the second form and applying M\"obius inversion:
\begin{equation}\label{gauss0}
I(z)=\sum_{k\ge1}\frac{\mu(k)}{k} \log\frac{1}{1-qz^k},
\end{equation}
and,  as a consequence, the number $I_n$ admits
the  explicit formula (already    known  to Gau{\ss},   see references
in~\cite[p.~46]{FlGoPa01}):
\begin{equation}\label{gauss}
I_n = \frac{1}{n}\sum_{k~|~n} \mu(k) q^{n/k} = \frac{q^n}{n}+R_n
\quad\hbox{where}\quad |R_n|\le q^{n/2}.
\end{equation}
In passing this approach also shows that 
\[
I(z)=\log\frac{1}{1-qz}+R(z),
\]
where   $R(z)$ is analytic  in  $|z|<q^{-1/2}$,  so  that $I(z)$  has a
logarithmic singularity at $z=q^{-1}$.

Let $D_n$ be the number of (monic) polynomials 
composed of irreducible factors of distinct degrees. The corresponding generating function, $D(z)$ 
satisfies
\begin{equation}\label{ddf1}
D(z)=\prod_{n=1}^\infty \left(1+I_n z^n\right),
\end{equation}
and the probability of the event over polynomials of degree~$n$ is
\[
\frac{D_n}{q^n} = [z^n]\, D\left(zq^{-1}\right).
\]
It is in fact the simplification  of this problem  for fields of large
cardinalities that originally led  Greene and Knuth   to consider the  problem of
enumerating  permutations     with  distinct         cycle     lengths
(see~\cite{FlGoPa01} for context). 
We have:

\begin{proposition} \label{poly-prop}
The probability that a random polynomial
of degree~$n$ over $\F_q$ has all its irreducible factors of distinct
degrees satisfies the asymptotic estimate
\begin{equation}\label{assert}
\frac{D_n}{q^n} = \delta(q) +O\left(\frac{1}{n}\right),
\end{equation}
where the constant $\delta(q)$ is given by
\begin{equation}\label{assert2}
\delta(q)= \prod_{k\ge1} \left(1+\frac{I_k}{q^k}\right) (1-q^{-k})^{I_k}.
\end{equation}
In addition, this probability admits a full asymptotic expansion with
oscillating coefficients.
\end{proposition}
 \begin{proof}
We start with a rough analysis based on comparing $D(z)$ and $P(z)$.
The second form of~\eqref{gau} implies
\begin{equation}
\frac{D(z/q)}{P(z/q)}=
\prod_{n\ge 1} \left(1+I_nq^{-n}z^n\right) (1-z^nq^{-n})^{I_n}.
\end{equation}
At $z=1$, the general factor in the product satisfies, for large~$n$,
\begin{equation}
\begin{array}{lll}\label{ddf3}
\ds (1+I_nq^{-n})(1-q^{-n})^{I_n}
& =& \ds \left(1+\frac1n+O(q^{-n/2})\right)
\exp\left(-\frac1n+O(q^{-n/2})\right)
\\
&=& \ds 1+\frac{1}{2n^2}+O(n^{-3}),
\end{array} 
\end{equation}
given the known form~\eqref{gauss} of~$I_n$. This is enough to
ensure the convergence of the infinite product defining $\delta(q)$
and proves that
\[
D\left(\frac{z}{q}\right)\sim \frac{\delta(q)}{1-z} \qquad (z\to1^-),
\]
which is evidently compatible with the asserted estimate~\eqref{assert}.
Also, the estimate~\eqref{ddf3} points to the fact that convergence 
of the infinite product giving~$\delta(q)$ is very slow, but also
suggests the modified scheme (compare with~\eqref{euler}),
\begin{equation}\label{ddf4}
\delta(q)= 
e^{-\gamma}
\prod_{k\ge1} \left(\frac{1+{I_k}{q^{-k}}}{1+1/k} 
(1-q^{-k})^{I_k} e^{1/k}\right). 
\end{equation}
which now exhibits geometric convergence.

In order to complete the derivation of~\eqref{assert}, we apply
the exp-log transformation 
to~\eqref{ddf1}:
\begin{equation}\label{ddf2}
D\left(\frac{z}{q}\right)=\exp\left(\sum_{m\ge1} \frac{(-1)^{m-1}}{m}\Lambda_m(z)\right),
\qquad
\Lambda_m(z):=\sum_{n\ge1} \left(I_nq^{-n}z^n\right)^m.
\end{equation}
Also, we have $\Lambda_1(z)=I(z/q)$. 
Taking out the factor $e^{\Lambda_1(z)}=e^{I(z/q)}$ gives
\begin{equation}\label{ddf5}\renewcommand{\arraycolsep}{2truept}
\begin{array}{lll}
\ds D(\frac{z}{q})&=&        \ds        e^{I(z/q)}\cdot         
\exp\left(\sum_{m\ge2} \frac{(-1)^{m-1}}{m}\Lambda_m(z)\right) \\
&=&\ds \frac{1}{1-z} \exp\left(\sum_{k\ge2} 
\frac{\mu(k)}{k}\log\frac{1}{1-q^{1-k}z^k}\right)
\exp\left(\sum_{m\ge2} \frac{(-1)^{m-1}}{m}\Lambda_m(z)\right)\\
&=&\ds \frac{e^{A(z)}}{1-z} \exp\left(\sum_{m\ge2} \frac{(-1)^{m-1}}{m}\Lambda_m(z)\right),
\end{array}
\end{equation}
where $A(z)$ is analytic in $|z|<q^{1/2}$.

Next, we observe from~\eqref{gauss} that 
\[
\Lambda_m(z)=\Li_m(z^m)+S_m(z),\]
where  $S_m(z)$ is   analytic in  $|z|<q^{1/2}$.   The  existence of a
factorization fulfilling  the conditions of Theorem~\ref{hybrid2-thm}
now   follows from the  usual split   of the   sum in the  exponential
of~\eqref{ddf5}, last line:
\[
\sum_{m\ge2} \frac{(-1)^{m-1}}{m} \Lambda_m(z)
= \sum_{m=2}^M \frac{(-1)^{m-1}}{m} \Lambda_m(z)
+ 
\sum_{m\ge M} \frac{(-1)^{m-1}}{m} \Lambda_m(z).
\]
The first split sum is singular at roots of unity of order $\le M$
and admits a log-power expansion of type $\cal O^t$ for any $t$;
the second split sum has its $n$th coefficient that is of order
$O(n^{-M-1})$, and is accordingly $\cal C^{M-1}$. 
By Theorem~\ref{hybrid2-thm}, these considerations imply that 
$[z^n]D(z/q)$ admits a full asymptotic expansion with
oscillating coefficients.

It is not hard then to track the order of the term that corrects
the dominant asymptotic regime: in essence, it arises from the factor
\[
\frac{1}{1-z} e^{-\frac12 \Li_2(z^2)},
\]
whose coefficients are $O(n^{-1})$.
\end{proof}

A variant of this problem is also considered in~\cite{FlGoPa01,KnWa95}:
\emph{determine the probability that a random polynomial of degree~$n$ is such that
its square-free factorization consists of polynomials of distinct degrees.} 
The generating function is then
\[
\widehat D(z)=\prod_{n\ge1} \left(1+\frac{I_nz^n}{1-z^n}\right).\]
By devices entirely similar to those employed in the proof of Proposition~\ref{poly-prop},
one finds that the corresponding probability satisfies
\[
\frac{\widehat{D}_n}{q^n}=\widehat\delta(q)+O\left(\frac{1}{n}\right),\]
where the constant $\widehat\delta(q)$ is 
\[
\widehat\delta(q)=
\prod_{k\ge1} \left(1+\frac{I_k}{q^k-1}\right) (1-q^{-k})^{I_k}
=e^{-\gamma}
\prod_{k\ge1} \left(\frac{1+{I_k}/{(q^{k}-1)}}{1+1/k} 
(1-q^{-k})^{I_k} e^{1/k}\right). 
\]
The treatment of~\cite{FlGoPa01} is based on a crude application
of the hybrid method, that of~\cite{KnWa95} relies on elementary bounds
and coefficient manipulations, but it is not clearly applicable to
derive asymptotic expansions beyond the first term.

\subsection{Dissimilar forests}
Fix a class $\cal T$ of rooted trees. An unordered forest 
is a multiset of trees in~$\cal T$, and we let $\cal F$ be the class of all
forests. A forest is said to be \emph{dissimilar} if
all the trees that it contains are of different sizes.
We let $\cal E$ denote the class of dissimilar forests.

For instance, 
the collection $\cal T$ of all plane trees
has generating function
\begin{equation}\label{for1}
T(z)=\frac{1-\sqrt{1-4z}}{2};
\end{equation}
a well-known fact~\cite{FlSe06,GoJa83,Knuth97,Stanley86,Wilf94}.
The generating function of unordered forests is
given by a construction analogous to~\eqref{gau}:
\[
\begin{array}{lll}
F(z)&=& \ds\prod_{n\ge1} (1-z^n)^{-T_n}
\quad=\quad \exp\left(\sum_{k=1}^\infty \frac{1}{k}T(z^k)\right)\\
&=&\ds 1+z+2z^2+4z^3+10z^4+26z^5+77z^6+235z^7+\cdots
\end{array}
\]
(the coefficients constitute~\EIS{A052854}).
The generating function of dissimilar forests, which is our actual object of study, is
\[
\begin{array}{lll}
E(z)&=& \ds \prod_{n=1}^\infty \left(1+T_n z^n\right)
\\
&=& 1+z+z^2+3z^3+7z^4+21z^5+63z^6+203z^7+\cdots
\end{array}\]
(currently unlisted in~\cite{Sloane06}).

The coefficients $[z^n]T(z)$ are the famous Catalan numbers of
combinatorial theory,
\[
T_n=\frac{1}{n}\binom{2n-2}{n-1}\sim \frac{4^{n-1}}{\sqrt{\pi n^3}}.\]
the asymptotic estimate being in agreement with the square-root singularity 
that is visible in~$T(z)$. An easy analysis of 
$F(z)$ then shows that, as $z\to\frac14$,
\begin{equation}\label{for3}
\begin{array}{lll}
F(z)& =& \ds K-\frac12 K \sqrt{1-4z}+\cdots\\
K&:=& \ds \exp\left(\sum_{k=1}^\infty \frac{1}{2k}
\left(1-\sqrt{1-4^{1-k}}\right)\right)\doteq 1.71603\,05349\,22281,
\end{array}
\end{equation}
resulting in the estimate 
\[
F_n\sim K \frac{4^{n-1}}{\sqrt{\pi n^3}}.
\]

\begin{proposition}\label{for-prop}
The probability that a random forest of rooted plane trees 
having total size~$n$ is dissimilar is
asymptotic to the constant
\begin{equation}\label{for4}
\kappa=\frac{e^{L+1/2}}{K} \qquad
\hbox{with}\quad L:=\sum_{m=2}^\infty \frac{(-1)^{m-1}}{m} \sum_{n\ge1}
\left(\frac{1}{n}\binom{2n-2}{n-1}4^{-n}\right)^m.
\end{equation}
This probability admits a complete asymptotic expansion with
oscillating coefficients.
\end{proposition}
\begin{proof}
We limit the calculation to first order asymptotics. 
As usual, start from the exp-log trick, which yields
\begin{equation}\label{for45}
E(z)=e^{T(z)} \exp\left(\sum_{m=2}^\infty \frac{(-1)^{m-1}}{m}\Lambda_m(z)\right),
\end{equation}
where the functions $\Lambda_m$ that play the r\^ole of special polylogarithms
are defined by
\begin{equation}\label{for5}
\Lambda_m(z):= \sum_{n\ge1} \left(\frac{1}{n}\binom{2n-2}{n-1}z^n\right)^m.
\end{equation}
Observe that the  functions $\Lambda_m(z)$  are $\Delta$-analytic,  by
virtue of a theorem asserting the  closure of functions of singularity
analysis  class  under Hadamard products~\cite{FiFlKa05,FlSe06}.  (The
function   $\Lambda_m(z)$  is  roughly  comparable  to $\Li_{3m/2}(4^m
z^m)$.)  Splitting  the sum in~\eqref{for45}  according to $m\le M$ and
$m>M$, for an arbitrarily large~$M$,  yields a factorization of $E(z)$
that fulfills the conditions  of Theorem~\ref{hybrid2-thm}.  Hence the
existence of a full expansion for coefficients is granted.

For first order  asymptotics, it suffices  to analyse  what happens as
$z\to(\frac14)^-$.  We find, from~\eqref{for45} and with the notations of~\eqref{for4}
\[
E(z)= e^{L+1/2} \left(1-\frac12\sqrt{1-4z}+O(1-4z)\right),
\]
which implies the stated estimate.
\end{proof}

The analysis obviously  extends to any simple family  of  trees, in the
sense of Meir and   Moon~\cite{FlSe06,MeMo78}, including the case  of  Cayley
trees that are enumerated by~$n^{n-1}$.

%
%
%

\section{Conclusion}\label{conc-sec}

As demonstrated by our foregoing examples, several infinite-product generating
functions occurring in combinatorics are amenable to the hybrid method,
despite the fact that they admit the unit circle as a natural boundary.
For comparison purposes, Figure~\ref{comp-fig} displays a few typical
infinite-product generating functions and the general asymptotic shapes\footnote{
We use `$\propto$' to represents asymptotic proportionality,
that is, asymptotic equivalence up to an unspecified multiplicative
constant.} of
their coefficients.

\begin{figure}\begin{small}
\frenchspacing
\begin{tabular}{l|l|l}
\hline
\emph{Generating function} & \emph{Type} & \emph{Coefficients} \\
\hline
\hline
$\ds Q_0=\prod_{n\ge1} \left(1+z^n\right)$ & Partitions dist. summ. & $\ds\propto
n^{-3/4}e^{\pi\sqrt{n/3}}$
\\
$\ds Q_{-1/2}=\prod_{n\ge1} \left(1+\frac{z^n}{\sqrt{n}}\right)$ & &
 $\propto
n^{-1/2}e^{\frac32\sqrt[3]{2\pi n}}$ \\
\hline
$\ds Q_{-1}=\prod_{n\ge1} \left(1+\frac{z^n}{n}\right)$ 
& Perm. dist. cycle lengths & $\ds\propto 1$
\\
$\ds Q_{-3/2}=\prod_{n\ge1} \left(1+\frac{z^n}{n^{3/2}}\right)$ 
& Forest dist. comp. sizes & $\ds\propto n^{-3/2}$
\\
$\ds Q_{-2}=\prod_{n\ge1} \left(1+\frac{z^n}{n^{2}}\right)$ 
& Perm. same cycle type & $\ds\propto n^{-2}$
\\
\hline
$\ds Q_{-\infty}=\prod_{n\ge1} \left(1+\frac{z^n}{n!}\right)$ 
& Set part. dist. block lengths & Mixed regimes~\cite{KnOdPiRi99}.\\
\hline
\end{tabular}
\end{small}
\caption{\label{comp-fig} Some typical infinite-product generating functions,
their ``types'',  and 
the asymptotic shape of the corresponding  coefficients.}
\end{figure}

Perhaps the most well-known generating function in this range is
\[
Q_0(z):=\prod_{n\ge1} \left(1+z^n\right),
\]
whose coefficients enumerate   partitions into distinct summands.   In
such a case, the analysis of coefficients, first performed by Hardy and
Ramanujan, is carried  out  by  mean  of the saddle  point  method, 
in
accordance with the fact that the  function $Q_0(z)$ gets 
exponentially large at its
``main''  singularity $z=1$.  This and similar  cases are best treated
by means of Meinardus' method~\cite{Andrews76}.

The function
\[
Q_{-1/2}(z):=\prod_{n\ge1} \left(1+\frac{z^n}{\sqrt{n}}\right),
\]
is also fast growing near $z=1$ and is best transformed by means of convergence factors:
\[
Q_{-1/2}(z)=\frac{e^{\Li_{1/2}(z)}}{\sqrt{1-z}}
\prod_{n\ge1} \left(1+\frac{z^n}{\sqrt{n}}\right)e^{-z/\sqrt{n}+z^2/2n}.\]
This expression shows a  growth  similar to that of the function
$\exp((1-z)^{-1/2})$, to which the saddle point method is applicable. 
All computations done, we find
\[
Q_{-1/2,n}\sim G\frac{e^{\zeta(1/2)}}{\sqrt{3\pi n}}
\exp\left(\frac32\sqrt[3]{2\pi n}\right)
\qquad\hbox{with}\quad 
G=\prod_{n\ge1}\left(1+\frac{1}{\sqrt{n}}\right)e^{-1/\sqrt{n}+1/(2n)}.\]

The next three examples of Figure~\ref{comp-fig} are relative to functions of the form
\[
Q_{\alpha}:=\prod_{n\ge1} \left(1+z^n n^{\alpha}\right),
\]
with $\alpha\le -1$. These are typical cases where the hybrid method is applicable.
The function $Q_{-1}$ is exactly the one encountered when enumerating 
permutations with distinct cycle lengths. The function
$Q_{-3/2}$ provides a simplified analytic model of dissimilar forests,
while $Q_{-2}$ behaves like the generating function of permutations 
having the same cycle type. The general pattern is, for $\alpha<-1$:
\[
Q_{\alpha,n}\propto n^{\alpha}.
\]

At the other end of the spectrum, we find the function 
\[
Q_{\infty}(z):=\prod_{n\ge1}\left(1+\frac{z^n}{n!}\right),
\]
which is the exponential generating function of set partitions, all of
whose  blocks are  of distinct   sizes. Clearly,   this is  an  entire
function to  which the hybrid  method is not  applicable. Coefficients
have been precisely studied by the seven authors of~\cite{KnOdPiRi99},
with some (but not all)
of the regimes being accessible to the saddle point method.

\smallskip

Finally, Gourdon has developed in~\cite{Gourdon96} methods leading
to asymptotic expansions of a shape similar to the ones of the present paper
in relation to a refinement of Golomb's
problem posed by Knuth in~\cite[Ex.~1.3.3.23]{Knuth97}.
He showed that the expected length of the longest cycle in a permutation of size~$n$
admits an asymptotic representation which starts as $(\omega=e^{2i\pi/3}$)
\[
\hbox{$ \lambda n+\frac12\lambda-\frac{e^\gamma}{24n}+
\big(\frac{e^\gamma}{48}-\frac{1}{8}(-1)^n\big)\frac{1}{n^2}
+\big(\frac{17e^\gamma}{3840}+\frac{1}{8}(-1)^n+\frac{1}{6}\omega^{1-n}+\frac{1}{6}\omega^{n-1}\big)\frac{1}{n^3}$,}
\]
with higher order roots of unity appearing at higher asymptotic orders. The proof
has analogies to the techniques of the present paper, with the  additional need of
precise expansions that describe the behaviour of the truncated 
logarithmic series on the unit disc.

\medskip

\noindent
{\bf Acknowledgements}. 
Bruno Salvy's {\tt equivalent} program has proved invaluable 
in permitting us to make several of our asymptotic expansions explicit.

\def\cprime{$'$}
\providecommand{\bysame}{\leavevmode\hbox to3em{\hrulefill}\thinspace}
\providecommand{\MR}{\relax\ifhmode\unskip\space\fi MR }
\providecommand{\MRhref}[2]{%
  \href{http://www.ams.org/mathscinet-getitem?mr=#1}{#2}
}
\providecommand{\href}[2]{#2}

\end{document}